\documentclass[11pt]{article}
\usepackage{amssymb,latexsym}
\usepackage{amsmath}

\title{ \large {\bf Finite basis problem for identities with involution.}}

\author{ Irina Sviridova\thanks{Supported by FAPESP, CNPq, CAPES;
e-mail \texttt{I.Sviridova@mat.unb.br}}
\\
\\
Departamento de Matem\'atica,\\
Universidade de Bras\'\i lia,\\
70910-900 Bras\'\i lia, DF, Brazil }

\date{October 26, 2014.}

\newtheorem{theorem}{Theorem}[section]
\newtheorem{definition}[theorem]{Definition}
\newtheorem{lemma}[theorem]{Lemma}

\newtheorem{remark}[theorem]{Remark}
\newtheorem{example}[theorem]{Example}
\newtheorem{conjecture}{Conjecture}[section]

\begin{document}
\maketitle

\begin{abstract}
We consider associative algebras
with involution over a field of characteristic zero.
We proved that any algebra with involution satisfies the same identities
with involution as the Grassmann envelope of some finite dimensional $(\mathbb{Z}/4 \mathbb{Z})$-graded algebra with graded involution. As a consequence we obtain the positive solution of the Specht problem for identities with involution: any associative algebra with involution over a field of characteristic zero has a finite basis of identities with involution. These results are analogs of Kemer's theorems for ordinary identities \cite{Kem1}. Similar results were proved also for associative algebras graded by a finite group in \cite{AB}, and for abelian case in \cite{Svi1}.
\bigskip


\end{abstract}

\section*{Introduction}

The interest to involutions on associative algebras can be partially explained
by their natural interconnections with various interesting and important classes
of algebras which appears in different fields of mathematics and physics
(see, e.g., \cite{Invbook}). Particularly, associative algebras with involution
is the natural background for important classes of Lie and Jordan algebras
(\cite{LieHum}, \cite{McCrim}, \cite{ZSSS}).
The identities with involution are also intensively studied last years.

In the theory of identities one of the central problem is the Specht problem.
This is the problem of existence of a finite base for any system of identities.
Originally this problem was formulated by W.Specht for ordinary polynomial
identities of associative algebras over a field of characteristic zero \cite{Specht}.
This problem was positively solved by A.Kemer \cite{Kem1}. The solution is based on the
Kemer's classification theorems. They state that any associative algebra over a field of
characteristic zero is equivalent in terms of identities (PI-equivalent) to the Grassmann envelope
of a finite dimensional superalgebra, and any finitely generated PI-algebra is PI-equivalent
to a finite dimensional algebra. The classification theorems have a proper significance.
They turn out the key tool for study of polynomial identities several last years.

The proof of the main classification theorem of Kemer consists of two principal steps:
the supertrick and the PI-representability of finitely generated PI-superalgebras.
On the first step the study of polynomial identities of any associative algebra is reduced
to study of identities of the Grassmann envelope of a finitely generated PI-superalgebra.
The second step is to prove that a finitely generated PI-superalgebra has the same $(\mathbb{Z}/2 \mathbb{Z})$-graded identities as some finite dimensional superalgebra.

Later results similar to some of the Kemer's theorems were obtained also for various classes
of algebras and identities. A review of results concerning the Specht problem can be found in
\cite{B6}. One of the most recent results is a positive solution of the local Specht problem
for associative algebras over an associative commutative Noetherian ring with unit \cite{B1}-\cite{B9}.
Graded algebras and algebras with involution were also considered with regard to this problem.
The positive solution of the Specht problem and analogs
of the classification theorems were obtained for graded identities of graded associative algebras
over a field of characteristic zero (\cite{AB} for a grading by a finite group,
and \cite{Svi1} for a grading by a finite abelian group). The equivalence
in terms of identities with involution was proved for finitely generated 
and finite dimensional PI-algebras with involution \cite{Svi2}.

The main purpose of this paper is a positive solution of the Specht problem for
identities with involution. This problem can be formulated in various forms:
in terms of a finite base of identities, and in terms of the Noetherian property for ideals of the free
algebra which are invariant under free algebra endomorphisms.
The positive answer to this question for identities with involution is equivalent to any of the following statement.
Any associative algebra with involution over a field of characteristic zero has a finite base
of identities with involution (all identities with involution of a $*$-algebra follow from a finite family of $*$-identities). Any $*$T-ideal of the free associative algebra with involution of infinite rank
over a field of characteristic zero is finitely generated as a $*$T-ideal. 
Any ascending chain of $*$T-ideals of the free associative algebra with involution of infinite rank
over a field of characteristic zero eventually stabilizes. 
$*$T-ideal is a $*$-invariant two-sided ideal of the free associative algebra with involution,
closed under all free algebra endomorphisms which commute with involution.
See Lemma \ref{*TS} about the structure of a $*$T-ideal generated by a set $S$.

We prove in this work that any associative algebra with involution over a field of characteristic zero satisfies the same identities with involution as the Grassmann $(\mathbb{Z}/4 \mathbb{Z})$-envelope of some finitely generated $(\mathbb{Z}/4 \mathbb{Z})$-graded PI-algebra with graded involution (Theorem \ref{IPI1}).
This is an analog of the supertrick in the classical case.
Using the recent result of the author about PI-representability of finitely generated $(\mathbb{Z}/4 \mathbb{Z})$-graded PI-algebras with graded involution \cite{Svi3} we obtain a version of the main
classification Kemer's theorem for identities with involution (Theorem \ref{IPI2}).
As a consequence we obtain the positive solution of the Specht problem for identities with involution
of associative $*$-algebras over a field of characteristic zero (Theorem \ref{IPI4}).

Throughout the paper we consider associative algebras over a field $F$ of characteristic zero.
Involution of an $F$-algebra $A$ is an anti-automorphism of $A$
of the second order. If we fix an involution $*$ of an associative
$F$-algebra $A$ then the pair $(A,*)$ is called an {\it associative
algebra with involution} (or {\it associative $*$-algebra}). Note that an
algebra with involution can be considered as an algebra with the
supplementary unary linear operation $*$ satisfying identities
\[
(a \cdot
b)^*=b^* \cdot a^*, \qquad (a^*)^*=a
\]
for all $a, b \in A.$

Observe that any $*$-algebra can be decomposed into
the sum of symmetric and skew-symmetric parts. An element $a \in
A$ is called {\it symmetric } if $a^*=a,$ and
{\it skew-symmetric } if $a^*=-a.$ So, $a+a^*$ is symmetric and
$a-a^*$ skew-symmetric for any $a \in A.$ Thus, we have
$A=A^{+} \oplus A^{-},$ where $A^{+}$ is the subspace
formed by all symmetric elements ({\it symmetric part}), and
$A^{-}$ is the subspace of all skew-symmetric elements of $A$
({\it skew-symmetric part}).
We also use the notations $a \circ b =a b + b a,$ and $[a,b]=a b - b a.$
It is clear that the symmetric part $A^{+}$ of
a $*$-algebra $A$ with the operation $\circ$ is a Jordan algebra (Hermitian Jordan algebra).
The skew-symmetric part $A^{-}$ with the operation $[,]$ is a Lie algebra.
All classical finite-dimensional simple Lie algebras over an algebraically closed
field, except $sl_n(F),$ are of this type \cite{LieHum}.

Suppose that $A,$ $B$ are algebras with involution. An ideal
$I$ of $A$ invariant with respect to involution is called {\it $*$-ideal.}
If $I \unlhd A$  is a $*$-ideal then $A/I$ inherits the involution of $A.$
A homomorphism $\varphi: A \rightarrow B$ is called {\it $*$-homomorphism }
({\it homomorphism of algebras with involution}) if it commutes with the involution.
We denote by $A_1 \times \dots \times A_{\rho}$ the direct product
of algebras $A_1, \dots, A_{\rho},$ and by $A_1 \oplus \dots
\oplus A_{\rho} \subseteq A$ the direct sum of subspaces $A_i$ of
an algebra $A.$ If $\tau_{i}$ is the involution
of $A_i$ ($i=1,\dots,\rho$) then $A_1 \times \dots \times A_{\rho}$ is an algebra with
the involution $*$ defined by the rules $(a_1,\dots,a_{\rho})^*=
(\tau_1(a_1),\dots,\tau_{\rho}(a_{\rho})),$ \ $a_i \in A_i.$

We study identities with involution ($*$-identities) of associative
algebras with involution. The notion of identity with involution is a formal extension of
the notion of ordinary polynomial identity (see, e.g., \cite{GZbook}, \cite{Svi2}).
A brief introduction to the notion is given in Section 1.
The definition of $*$-identity can be found also in \cite{Svi2} or in \cite{GZbook}
with some more details.
We refer the reader to the textbooks \cite{Dren},
\cite{DrenForm}, \cite{GZbook}, and to \cite{BelRow}, \cite{Kem1} concerning basic
definitions, facts and properties of ordinary polynomial identities.

We also use in the proof of the classification theorem the concept of a graded identity
with involution (graded $*$-identity).
This concept was developed in \cite{Svi3}. The principal definitions concerning this notion
is also given in Section 1. In general, the concept of a graded $*$-identity is the union
of concepts of an identity with involution and of a graded identity. The information about graded
identities can be found in \cite{GRZ1}, \cite{GZbook} and in \cite{AB}, \cite{Svi1}.

Besides the notions of the free algebra with involution, identities with involution,
the free graded algebra with involution and graded identities with involution, Section 1 also
contains the necessary information about graded algebras.

Properties of multilinear $*$-polynomials and multilinear graded $*$-polynomials alternating or
symmetrizing in some set of variables are discussed in Section 3. Such polynomials appears in
the study of identities as a result of application of techniques of symmetric group representations.
Basic facts and notions concerning applications of representation theory for $*$-identities
can be found in \cite{DrenGiam}, \cite{GR}, \cite{GM}, \cite{GM1}, \cite{GZbook}.
Observe that in our case the application of representation theory for $*$-identities
is similar to the case of ordinary polynomial identities due to fact
that the symmetric group acts by renaming of variables on a homogeneous subset
of variables (on a set of symmetric variables in respect to involution or skew-symmetric).
Thus in many situations we can apply the same results and arguments as in the case
of ordinary polynomial identities.
The book \cite{GZbook} contains very detailed and complete exposition of the facts and methods related to application of symmetric group representations for theory of polynomial identities.
We appeal to this book when we need facts which can be directly applied in our case
or arguments which can be literally repeated.
We also refer the reader to \cite{JKer}, \cite{CurtR} concerning principal definitions
and facts of representation theory.

Section 2 is devoted to the definition of the Grassmann envelope of a
$(\mathbb{Z}/4 \mathbb{Z})$-graded algebra. Section 4 contains the classification theorems
for ideals of identities with involution (Theorems \ref{IPI1}, \ref{IPI2}, \ref{IPI3}).
They are analogs of Kemer's theorems \cite{Kem1}
for polynomial identities of associative algebras over a field of characteristic zero.
The proof of Theorem \ref{IPI1} follow the scheme of the proof
of the classical Kemer's theorem about Grassmann envelopes given in \cite{GZbook}.
We adopt this proof for the case of identities with involution. Theorem \ref{IPI2} is the
corollary of Theorem \ref{IPI1} and Theorem 6.2 \cite{Svi3}.
The Specht problem solution (Theorem \ref{IPI4}) for $*$-identities is given in Section 5.
The proof of Theorem \ref{IPI4} is the involution version of the original
Kemer's proof \cite{Kem1}.

Observe that the principal tool of the proof is the Grassmann envelope.
Our conception of the Grassmann envelope in this work is different of the usual one.
Usually one consider the Grassmann envelope $E(A)=A_{\bar{0}} \otimes E_{\bar{0}} \oplus A_{\bar{1}} \otimes E_{\bar{1}}$ for a $(\mathbb{Z}/2 \mathbb{Z})$-graded algebra $A=A_{\bar{0}} \oplus A_{\bar{1}}$ (superalgebra). It gives super-theory. In this case a graded involution on $E(A)$ induces the superinvolution
on $A.$ A $(\mathbb{Z}/2 \mathbb{Z})$-graded linear transformation $\star$ of the second order of a superalgebra $A$ is called a superinvolution if
\[ (a \cdot b)^{\star} = (-1)^{i \cdot j} \
b^{\star} a^{\star} \qquad \forall a \in A_{\bar{i}}, \  b \in A_{\bar{j}}, \ \ i, j \in \{ 0, 1\}.\]
And vice versa, one needs a superinvolution on $A$ to guarantee the
correspondent involution on $E(A).$

We use a slight generalization of the traditional construction based on the natural
$(\mathbb{Z}/4 \mathbb{Z})$-grading of the Grassmann algebra $E.$ We call it
the Grassmann $\mathbb{Z}_4$-envelope to differ it from the traditional Grassmann envelope.
This construction is compatible with the usual graded involution. We think that the Specht problem for
$*$-identities can be solved also using the traditional approach based on superinvolutions.
It is possible even that the traditional approach could be more natural.
But the author assume that the new construction and its connection with graded involutions
on associative algebras is rather curious and worth to study.

\section{Identities with involution and graded identities with involution.}

Let $F$ be a field of characteristic zero.
Consider two countable sets $Y = \{ y_{i} | i \in \mathbb{N} \},$ \ $Z = \{ z_{i}
| i \in \mathbb{N} \}$
of pairwise different letters, and the free associative non-unitary
algebra $F\langle Y, Z \rangle$ generated by $Y \cup Z.$
We can define an involution on $F\langle Y, Z \rangle$ assuming that
variables from $Y$ are symmetric, and from $Z$ skew-symmetric
\begin{eqnarray} \label{freeinv}
&&(\sum \alpha_{w} \  a_{i_1} \cdots a_{i_n})^*=\sum \alpha_{w}\   a_{i_n}^* \cdots a_{i_1}^*
=\sum (-1)^{\deg_{Z} w}  \alpha_{w} \  a_{i_n} \cdots a_{i_1}, \quad \mbox{where} \nonumber\\
&&y_{j}^*=y_{j}, \ \ z_{j}^*=-z_{j},
\qquad w=a_{i_1} \cdots a_{i_n}, \ a_j \in Y \cup Z, \ \alpha_{w} \in F.
\end{eqnarray}
$F\langle Y, Z \rangle$ is the free associative algebra with
involution. Its elements are called {\it $*$-polynomials}.
The free associative algebra $F \langle X^* \rangle$ generated by the set
$X^*=\{ x_i, x_i^* | i \in \mathbb{N} \}$ also has an involution defined by
\begin{eqnarray*}
&&(\sum \alpha_{w} \  a_{i_1} \cdots a_{i_n})^*=\sum \alpha_{w} \  a_{i_n}^* \cdots a_{i_1}^*, \quad \mbox{where} \nonumber\\
&&(x_{j})^*=x_{j}^*, \ \ (x_{j}^*)^*=x_{j},
\qquad w=a_{i_1} \cdots a_{i_n}, \  a_j \in X^{*}, \  \alpha_{w} \in F.
\end{eqnarray*}
The equalities
\begin{eqnarray} \label{free}
&&y_i=\frac{x_i+x_i^*}{2}, \ \  z_i=\frac{x_i-x_i^*}{2}; \nonumber
\\ &&x_i=y_i+z_i, \ \ x_i^*=y_i-z_i
\end{eqnarray}
induce the isomorphism of algebras with involution $F\langle
X^*\rangle$ and $F\langle Y, Z \rangle.$ We use the algebra $F\langle Y, Z \rangle$
as the free associative $*$-algebra.

An algebra with involution $A$
satisfies {\it the $*$-identity } (or {\it identity with involution})
$f=0$ for a non-trivial $*$-polynomial
$f=f(y_{1},\dots,y_{n},z_{1},\dots,z_{m}) \in F\langle Y, Z \rangle$
whenever $f(a_1,\dots,a_n,b_1,\dots,b_m)=0$ for all elements
$a_i \in A^{+},$ and $b_i \in A^{-}.$
Let $\mathrm{Id}^{*}(A)$ be the ideal of all identities with involution of $A.$
Then $\mathrm{Id}^{*}(A)$ is a two-sided $*$-ideal of $F\langle Y, Z \rangle$ closed under all 
$*$-endomorphisms of $F\langle Y, Z \rangle$. Such ideals are called {\it $*$T-ideals} (see \cite{Svi2}).
Conversely, any $*$T-ideal $I$ of $F\langle Y, Z \rangle$ is the ideal of
$*$-identities of the algebra with involution $F\langle Y, Z \rangle/I.$
We denote by $*T[S]$ the $*$T-ideal generated by a set $S \subseteq F\langle Y, Z \rangle.$
The next statement is clear due to the definition and elementary properties of a $*$T-ideal.

\begin{lemma} \label{*TS}
Let $F$ be a field of characteristic zero.
Given a set $S \subseteq F\langle Y, Z \rangle$
a polynomial $f \in F\langle Y, Z \rangle$ belongs
to the $*$T-ideal $*T[S]$ generated by $S$ iff $f$ is a finite linear
combination of the form
\begin{equation}\label{form-gen}
f=\sum_{(u), j} \alpha_{(u), j}
\  v_{1} g_{j}(\tilde{u}_{j 1},\dots,\tilde{u}_{j n_{j}}) v_{2}, \ \alpha_{(u), j} \in F.
\end{equation}
Where $g_j=\tilde{g}_j$ or
$g_j=\tilde{g}^{*}_j$ for the full linearization $\tilde{g}_j$ of a multihomogeneous component
of a polynomial $g \in S;$ \ $\tilde{u}_{j l} = u_{j l} \pm u_{j l}^*$ for a
monomial $u_{j l} \in F\langle Y, Z \rangle$ ($\tilde{u}_{j l} = u_{j l} + u_{j l}^*$
if the corresponding variable $x_{j l}$ of the polynomial $g_j$ is symmetric
in respect to involution ($x_{j l} \in Y$), and $\tilde{u}_{j l} = u_{j l} - u_{j l}^*$
if $x_{j l} \in Z$ is skew-symmetric); \
$v_{l} \in F\langle Y, Z \rangle$ are
monomials, possibly empty; \ $(u)=(v_{1}, \tilde{u}_{j 1},\dots,\tilde{u}_{j n_{j}},v_{2}).$
\end{lemma}
\noindent {\bf Proof.}
It is clear that the set of all polynomials of the form (\ref{form-gen})
is a $*$T-ideal, and contains $S.$ The characteristic of the base field is zero.
Therefore, any $*$T-ideal $\Gamma$ contains all multihomogeneous components
of its elements and their full linearizations. Particularly, if a $*$T-ideal $\Gamma$
contains $S$ then it contains also all multihomogeneous components
of any $g \in S$ and their full linearizations. Moreover, any $*$-invariant evaluation
of variables of a homogeneous polynomial $\widetilde{g} \in \Gamma$ can be realized by a $*$-invariant evaluation of the full linearization of $\widetilde{g}$ up to a non-zero coefficient.
Since the polynomials $g_{i}$ are multilinear then a linear base of all their $*$-invariant evaluations
is formed by their evaluations with the symmetric and skew-symmetric parts of monomials.
\hfill $\Box$

We say that a $*$-polynomial $f$ is a {\it consequence of a set $S \subseteq F\langle Y, Z \rangle$} if 
$f \in *T[S].$ 
We have also that $\mathrm{Id}^{*}(A_1 \times A_{2}) = \mathrm{Id}^{*}(A_1) \cap 
\mathrm{Id}^{*}(A_2)$  for the direct product $A_1 \times A_{2}$ of arbitrary
$*$-algebras $A_i.$

Suppose that $\Gamma$ is a $*$T-ideal. {\it A $*$-variety defined by $\Gamma$}
is the family of all associative $*$-algebras such that they satisfy $f=0$ for any $f \in \Gamma.$
It is denoted by $\mathfrak{V}_{\Gamma}.$
A $*$-algebra $A$ generates $\mathfrak{V}_{\Gamma}$ if $\Gamma=\mathrm{Id}^{*}(A).$ Then we write
$\mathfrak{V}_{\Gamma}=\mathfrak{V}(A).$
The $*$-algebra $F\langle Y, Z \rangle / \Gamma$ is the relatively free
algebra of the $*$-variety $\mathfrak{V}_{\Gamma}.$
Any $*$-variety is closed under taking $*$-subalgebras,
$*$-homomorphic images, and direct products.
The free $*$-algebra of rank $\nu$ $F\langle Y_{\nu}, Z_{\nu} \rangle,$
and the relatively free algebra of rank $\nu$
$F\langle Y_{\nu}, Z_{\nu} \rangle /( \Gamma \cap F\langle Y_{\nu}, Z_{\nu} \rangle)$
for the $*$-variety $\mathfrak{V}_{\Gamma}$ are also considered 
($Y_{\nu} = \{ y_{i} | i=1,\dots,\nu \},$ \ $Z_{\nu} = \{ z_{i}
| i=1,\dots,\nu \}$).

Let $G$ be a finite abelian group.
An algebra $A$ is {\it $G$-graded } if $A=\bigoplus_{\theta \in G} A_{\theta}$
is the direct sum of its subspaces $A_{\theta}$ satisfying
$A_{\theta} A_{\xi} \subseteq A_{\theta \xi}$ for all $\theta, \xi \in G.$
An element $a \in A_{\theta}$ is called {\it $G$-homogeneous of degree }
$\deg_{G} a=\theta.$ A subspace $V$ of $A$ is graded
if $V=\bigoplus_{\theta \in G} (V \cap A_\theta).$

\begin{example} \label{ZnGR}
The free associative algebra $\mathfrak{F}=F\langle X \rangle$ generated by
$X=\{ x_1, x_2, \dots \}$ has the natural $(\mathbb{Z}/n \mathbb{Z})$-grading
$\mathfrak{F}_{\bar{m}}=\mathrm{Span}_F\{
x_{i_1} x_{i_2} \cdots x_{i_s} | s \equiv m \ \mathrm{mod } \ n \},$  \
$\bar{m} \in \mathbb{Z}/n \mathbb{Z}.$

The Grassmann algebra of countable rank \  $E=\langle
e_i, \  i \in \mathbb{N} |\  e_i e_j = -e_j e_i, \ \forall i, j
\rangle$ has the homogeneous relations. Thus it inherits the $(\mathbb{Z}/n \mathbb{Z})$-grading
of the free algebra $E_{\bar{m}}=\mathrm{Span}_F\{
e_{i_1} e_{i_2} \cdots e_{i_s} | s \equiv m \ \mathrm{mod } \ n, \  i_1 < \dots < i_s \}.$
This grading is called {\it natural}.
\end{example}

Consider a $G$-graded algebra $A$ with involution.
We assume that the involution is a graded anti-automorphism of $A,$ i.e.
$A_{\theta}^{*}=A_{\theta}$ for any $\theta \in G.$ This is equivalent to condition
(see for instance \cite{BahtShestZ}) that the subspaces
$A^{+},$ $A^{-}$ are graded. Particularly, we have that
$A=\bigoplus_{\theta \in G} (A_{\theta}^{+} \oplus A_{\theta}^{-}),$ where
$A^{\delta}=\bigoplus_{\theta \in G} A_{\theta}^{\delta},$
\ ($\delta \in \{ +, -\}$); and
$A_{\theta}=A_{\theta}^{+} \oplus A_{\theta}^{-},$ \ ($\theta \in G$).
We say that an element $a \in A_{\theta}^{\delta}$ ($\delta \in \{ +, -\},$ \
$\theta \in G$) is {\it homogeneous of complete degree }
$\deg_{\widehat{G}} a=(\delta,\theta)$ or simply {\it $\widehat{G}$-homogeneous}.

\begin{example} \label{GrassmInv}
Consider the natural $(\mathbb{Z}/4 \mathbb{Z})$-grading on
the Grassmann algebra of countable rank $E=\bigoplus_{\bar{m} \in \mathbb{Z}/4 \mathbb{Z}} E_{\bar{m}}$
described in Example \ref{ZnGR}. Define on $E$ the involution $*_E$ by the equalities
$(e_i)^{*_E}=e_i$ for all $i \in \mathbb{N}$.
This involution is called {\it canonical}. It is clear that this involution is graded.
Moreover, $E^{+}=E_{\bar{0}} \oplus E_{\bar{1}},$ and $E^{-}=E_{\bar{2}} \oplus E_{\bar{3}}.$
\end{example}

A homomorphisms $\varphi: A \rightarrow B$ of two $G$-graded $*$-algebras
$A,$ $B$ is called {\it graded $*$-homomorphism } if $\varphi$ is graded
($\varphi(A_{\theta}) \subseteq A_{\theta}$ \ for any $\theta \in G$),
and commutes with the involution. An ideal (a subalgebra) $I \unlhd A$
of a graded algebra with involution $A$ is {\it graded $*$-ideal }
({\it graded $*$-subalgebra }) if it is graded and invariant under the involution.
For graded algebras with involution we
consider only graded $*$-ideals, and graded $*$-homomorphisms.
In this case the quotient algebra $A/I$ is also a graded $*$-algebra
with the grading and the involution induced from $A$.
It is clear that the direct product of graded algebras
with involution is also a graded algebra with involution
(the grading and the involution are component-wise).

We can also define the notion of a graded $*$-identity for a $G$-graded algebra with
a graded involution.
The free associative non-unitary algebra $\mathfrak{F}^{G}=F\langle Y^{G}, Z^{G} \rangle$ generated by the set
$Y^{G} \cup Z^{G} = \{ y_{i \theta} | \theta \in G, i \in \mathbb{N} \} \cup \{ z_{i \theta}
| \theta \in G, i \in \mathbb{N} \}$ has the involution defined by
(\ref{freeinv}) for monomials in $Y^{G} \cup Z^{G}.$ We assume that
$y_{j \theta}^*=y_{j \theta},$ and $z_{j \theta}^*=-z_{j \theta}$
(for all $\theta \in G,$ $i \in \mathbb{N}$).
The $G$-grading on $\mathfrak{F}^{G}$ is defined naturally by the rule
\ $\deg_{G} a_{i_1} a_{i_2} \cdots a_{i_n} =
\deg_{G} a_{i_1} \cdots \deg_{G} a_{i_n},$ \ where $\deg_{G} y_{i \theta} = \deg_{G} z_{j \theta} =
\theta,$ \ $a_{j} \in Y^{G} \cup Z^{G}.$
It is clear that the involution (\ref{freeinv}) is graded.
The algebra $\mathfrak{F}^{G}$ is the free associative $G$-graded algebra with
graded involution. Its elements are called {\it graded $*$-polynomials}.
Variables $y_{i \theta} \in Y^{G},$ $z_{j \theta} \in Z^{G}$ are
$\widehat{G}$-homogeneous. Their complete degrees are $\deg_{\widehat{G}} y_{i \theta}=(+,\theta),$
\ $\deg_{\widehat{G}} z_{i \theta}=(-,\theta),$ \ $\theta \in G.$
Let us denote also $Y_{\theta}=\{ y_{i \theta} | i \in \mathbb{N} \},$ and
$Z_{\theta}=\{ z_{i \theta} | i \in \mathbb{N} \}$ for any $\theta \in G.$

Let $f=f(x_{1},\dots,x_{n}) \in F\langle Y^{G}, Z^{G} \rangle$ be a
non-trivial graded $*$-polynomial ($x_i \in Y^{G} \cup Z^{G}$).
We say that a graded $*$-algebra $A$
satisfies {\it the graded $*$-identity} (or {\it graded identity with involution})
$f=0$ iff $f(a_1,\dots,a_n)=0$ for all $\widehat{G}$-homogeneous elements
$a_i \in A_{\theta_i}^{\delta_i}$ of the corresponding complete degrees
$\deg_{\widehat{G}} a_i=\deg_{\widehat{G}} x_i=
(\delta_i,\theta_i),$ \ $\delta_i \in \{ +, - \},$ \
$\theta_i \in G$ ($i=1,\dots,n$).

Denote by $\mathrm{Id}^{gi}(A) \unlhd F\langle Y^{G}, Z^{G} \rangle$  the
ideal of all graded identities with involution of a graded $*$-algebra $A.$
It is clear that $\mathrm{Id}^{gi}(A)$
is a two-side graded $*$-ideal of $F\langle Y^{G}, Z^{G} \rangle$ closed under 
graded $*$-endomorphisms of $F\langle Y^{G}, Z^{G} \rangle$. 
We call such ideals {\it $gi$T-ideals} (see \cite{Svi3}). Conversely, any
$gi$T-ideal $I$ of $F\langle Y^{G}, Z^{G} \rangle$ is the ideal of graded
$*$-identities of the graded algebra with involution $F\langle Y^{G}, Z^{G} \rangle/I.$
Given a set $S \subseteq F\langle Y^{G}, Z^{G} \rangle$ of graded $*$-polynomials
denote by $giT[S]$ the
$gi$T-ideal generated by $S.$ Similarly to case of non-graded $*$-identities,
we have that $\mathrm{Id}^{gi}(A_1 \times \dots \times A_{\rho}) = \bigcap \limits_{i=1}^{\rho}
\mathrm{Id}^{gi}(A_i)$  for the direct product $A_1 \times \dots \times A_{\rho}$ of graded
$*$-algebras.

Given a $gi$T-ideal $\Gamma$ consider the family $\mathfrak{V}^{G}_{\Gamma}$
of all associative $G$-graded $*$-algebras that satisfy $f=0$ for any $f \in \Gamma.$
We call $\mathfrak{V}^{G}_{\Gamma}$ {\it a graded $*$-variety defined by $\Gamma.$}
If $\Gamma=\mathrm{Id}^{gi}(A)$ then we say that the graded $*$-algebra $A$
generates the graded $*$-variety $\mathfrak{V}^{G}_{\Gamma}=\mathfrak{V}^{G}(A).$
Particularly, $\mathfrak{V}^{G}_{\Gamma}=\mathfrak{V}^{G}(F\langle Y^{G}, Z^{G} \rangle / \Gamma).$
Moreover, the algebra $\mathfrak{F}_{\Gamma}=F\langle Y^{G}, Z^{G} \rangle / \Gamma$ is the relatively free
algebra of the graded $*$-variety $\mathfrak{V}^{G}_{\Gamma}.$
It is clear that $B \in \mathfrak{V}^{G}(A)$ for a graded $*$-algebra $B$ whenever
$\mathrm{Id}^{gi}(A) \subseteq \mathrm{Id}^{gi}(B).$
Any graded $*$-variety is closed under taking graded $*$-subalgebras,
graded $*$-homomorphic images, and direct products.

Let $Y^{G}_\nu = \{ y_{i \theta} | \theta \in G, 1 \le i \le \nu \},$
$Z^{G}_\nu = \{ y_{i \theta} | \theta \in G, 1 \le i \le \nu \}$ be two finite sets,
$\nu \in \mathbb{N}.$
We also consider the free $G$-graded algebra with involution
$F\langle Y^{G}_{\nu}, Z^{G}_{\nu} \rangle$ of rank $\nu$
generated by $Y^{G}_\nu \cup Z^{G}_\nu$ and the
relatively free algebra of rank $\nu$ $\mathfrak{F}_{\nu,\Gamma}=
F\langle Y^{G}_{\nu}, Z^{G}_{\nu} \rangle /( \Gamma \cap F\langle Y^{G}_{\nu}, Z^{G}_{\nu} \rangle)$
for the graded $*$-variety $\mathfrak{V}^{G}_{\Gamma}.$

Observe that omitting indices by the elements of the group $G$ in the structures
of the free graded $*$-algebra, graded $*$-identities and graded $*$-varieties
we obtain the notions of non-graded identities with involution and non-graded $*$-varieties.
Notice that in both cases (graded and non-graded) variables of the set $Y$ are reserved
for symmetric elements, and variables $Z$ for skew-symmetric.
Two $G$-graded algebras with involution $A$ and $B$ are called
{\it $gi$-equivalent}, $A \sim_{gi} B,$ if $\mathrm{Id}^{gi}(A)=\mathrm{Id}^{gi}(A).$
Non-graded algebras with involution $A$ and $B$ are
{\it $*$PI-equivalent}, $A \sim_{*} B,$ if
$\mathrm{Id}^*(A)=\mathrm{Id}^*(B).$
Given a $gi$T-ideal ($*$T-ideal) $\Gamma$ and
graded (non-graded) $*$-polynomials $f, g$ we write $f=g \
(\mathrm{mod } \ \Gamma)$ if $f-g \in
\Gamma.$

If we have a graded $*$-algebra $A$ then we assume
that $\mathrm{Id}^{*}(A) \subseteq \mathrm{Id}^{gi}(A).$
Namely, for a non-graded $*$-polynomial $f(y_{1},\dots,y_{n},z_{1},\dots,z_{m}) \in F\langle Y, Z \rangle$
we assume $f \in \mathrm{Id}^{gi}(A)$ whenever
$f(\sum_{\theta \in G} y_{1 \theta},\dots,\sum_{\theta \in G} y_{n \theta}, \sum_{\theta \in G} z_{1 \theta},\dots,\sum_{\theta \in G} z_{m \theta})
\in \mathrm{Id}^{gi}(A).$ Particularly, for a multilinear non-graded $*$-polynomial
$f(y_{1},\dots,y_{n},z_{1},\dots,z_{m}) \in F\langle Y, Z \rangle$
we have $f \in \mathrm{Id}^{gi}(A)$ if and only if
$f(y_{1 \theta_1},\dots,y_{n \theta_n},z_{1 \theta_{n+1}},\dots,z_{m \theta_{n+m}})
\in \mathrm{Id}^{gi}(A)$ for all $(\theta_1,\dots,\theta_{n+m}) \in G^{n+m}.$
Thus if $A \sim_{gi} B$ for $G$-graded $*$-algebras $A,$ $B$ then we have also
that $A \sim_{*} B.$

Note that the set $X^{G} = \{ x_{i \theta}=y_{i \theta}+z_{i \theta}
| i \in \mathbb{N}, \  \theta \in G \}$ generates in $\mathfrak{F}^{G}$
a $G$-graded subalgebra $F\langle X^{G} \rangle$ which is 
isomorphic to the free associative $G$-graded algebra
(\cite{Svi1}). Thus the ideal $\mathrm{Id}^{G}(A)$ of graded
identities of $A$ also lies in $\mathrm{Id}^{gi}(A).$

Recall that an algebra $A$ is called {\it PI-algebra } if it satisfies
a non-trivial ordinary polynomial identity (non-graded and without involution)
(see \cite{Dren}, \cite{DrenForm}, \cite{GZbook}, \cite{BelRow}, \cite{Kem1}).
It is clear that for a $G$-graded PI-algebra $A$ with involution
the T-ideal of ordinary polynomial identities $\mathrm{Id}(A)$ also lies in $\mathrm{Id}^{gi}(A).$
Moreover, we have that $\mathrm{Id}(A) \subseteq \mathrm{Id}^{*}(A) \subseteq \mathrm{Id}^{gi}(A).$
Here for a polynomial $f(x_1,\dots,x_n) \in \mathrm{Id}(A)$ we assume that
$f \in \mathrm{Id}^*(A)$ iff $f(y_1+z_1,\dots,y_n+z_n) \in \mathrm{Id}^*(A).$
This is the natural relation induced by the isomorphism \eqref{free} of $F\langle
X^*\rangle$ and $F\langle Y, Z \rangle$ and inclusion $F\langle
X \rangle \subseteq F\langle X^*\rangle.$

By Amitsur's theorem \cite{A3}, \cite{A4} (see also \cite{GZbook})
any $*$-algebra satisfying a non-trivial $*$-identity is a PI-algebra.
Thus any non-trivial $*$T-ideal contains a non-trivial T-ideal.
A $G$-graded $*$-algebra can not be a PI-algebra in general
(see for instance comments after Theorem 1 \cite{Svi1}).
In general case a graded $*$-algebra $A$ is a PI-algebra iff the neutral component
$A_{\mathfrak{e}}$ satisfies a non-trivial $*$-identity, where $\mathfrak{e}$ is the unit element of $G$
(it follows from \cite{A3}, \cite{A4}, and \cite{BahtGR}, \cite{BergC}).
This is equivalent to condition that $A$ satisfies a non-trivial non-graded $*$-identity.

The notion of degree of a graded or non-graded $*$-polynomial
is defined in the usual way. Using the multilinearization process
as in the case of ordinary identities (\cite{Dren}, \cite{DrenForm}, \cite{GZbook})
we can show that any $gi$T-ideal or $*$T-ideal over a field of characteristic
zero is generated by multilinear polynomials (see also Lemma \ref{*TS}).
Thus in our case it is enough to consider only multilinear identities.

The space of multilinear $*$-polynomials of degree $n$ has
the form
\[P_n=\mathrm{Span}_F\{ x_{\sigma(1)} \cdots x_{\sigma(n)}
| \sigma \in \mathrm{S}_n, x_i \in Y \cup Z \}.\]
Thus $P_n$ is the direct sum of subspaces of multihomogeneous and
multilinear polynomials depending on a fixed set of symmetric and
skew-symmetric variables. When we consider $*$-identities
we can assume that a multilinear $*$-identity depends on variables
$\{y_{1}, \dots, y_{k}\},$ and $\{ z_{1}, \dots, z_{n-k} \},$ \
$k=0,\dots,n.$ Denote by $P_{k,n-k}$ the subspace of all multilinear
$*$-polynomials $f(y_1,\dots,y_{k},z_1,\dots,z_{n-k})$ for a fixed number $k.$
Given a $*$T-ideal $\Gamma \unlhd F\langle Y, Z \rangle$
the vector spaces $\Gamma_{k,n-k}=\Gamma \cap P_{k,n-k},$
and $P_{k,n-k}(\Gamma)=P_{k,n-k}/
\Gamma_{k,n-k} \subseteq F\langle Y, Z \rangle/\Gamma$
has the natural structure of $(FS_{k} \otimes FS_{n-k})$-modules.
Here $S_{k}$ and $S_{n-k}$ act on symmetric and
skew-symmetric variables independently renaming the variables (see, e.g., \cite{GM}).

Further we consider $(\mathbb{Z}/4 \mathbb{Z})$-graded algebras with involution
and $(\mathbb{Z}/4 \mathbb{Z})$-graded $*$-identities. We assume that $G=\mathbb{Z}/4 \mathbb{Z},$
and use for it the additive notation. 
We also denote for brevity the group $\mathbb{Z}/4 \mathbb{Z}$ by $\mathbb{Z}_4,$ and
the free $\mathbb{Z}_4$-graded $*$-algebra $F\langle Y^{\mathbb{Z}_4}, Z^{\mathbb{Z}_4} \rangle$ by
$\mathfrak{F}^{(4)}.$ 

Let us define the function
$\eta: \mathbb{Z}_4 \rightarrow \{0, 1 \}$ 
by the rules $\eta(\bar{0})=\eta(\bar{1})=0,$
$\eta(\bar{2})=\eta(\bar{3})=1.$ The next elementary properties of $\eta$ can be checked directly
\begin{eqnarray*}
&&\eta(x)+\eta(y)=\eta(x+y) + 1 \ {\rm mod } \ 2 \quad \ \mbox{ if } \  x, y \in \{ \bar{1}, \bar{3}\}, \\
&&\eta(x)+\eta(y)=\eta(x+y) \ {\rm mod } \  2 \quad \ \mbox{ if } \  x \ \mbox{ or } \  y \  \mbox{ is even.}
\end{eqnarray*}

\section{Grassmann $\mathbb{Z}_4$-envelope of a graded $*$-algebra.}

Assume that $G=\mathbb{Z}_4.$ Consider a $\mathbb{Z}_4$-graded algebra
$A=\bigoplus_{\theta \in \mathbb{Z}_4} A_{\theta}.$

\begin{definition}
The algebra $E_4(A)=\bigoplus_{\theta \in \mathbb{Z}_4} A_{\theta} \otimes_F E_{\theta}$
is called Grassmann $\mathbb{Z}_4$-envelope of $A.$
Where $E=\bigoplus_{\theta \in \mathbb{Z}_4} E_{\theta}$ is the natural
$\mathbb{Z}_4$-grading of $E$ defined in Example \ref{ZnGR}.
\end{definition}
The algebra $E_4(A)$ is also $\mathbb{Z}_4$-graded with the grading
$(E_4(A))_{\theta}=A_{\theta} \otimes_F E_{\theta},$ \ $\theta \in \mathbb{Z}_4.$
If $A$ has a graded involution $*_A$ then the $F$-linear involution $*$ on
$E_4(A)$ is defined by the rules $(a \otimes g)^*=a^{*_A} \otimes g^{*_E},$
where $*_E$ is the canonic involution on $E$ (see Example \ref{GrassmInv}).
Hence $(a_{\theta} \otimes g_{\theta})^*=(-1)^{\eta(\theta)} \  a_{\theta}^{*_A} \otimes g_{\theta}$
for any $a_{\theta} \in A_{\theta},$ \ $g_{\theta} \in E_{\theta},$ \ $\theta \in \mathbb{Z}_4.$
It is clear that $E_4(A)^{\delta}=\bigoplus_{\theta \in G} (E_4(A))_{\theta}^{\delta},$ \ $\delta \in \{ +, -\},$
where
\begin{eqnarray} \label{Esym}
&&(E_4(A))_{\theta}^{+}=\mathrm{Span}_F\{
a_{\theta} \otimes g_{\theta} | a_{\theta} \in A_{\theta}^{+}, \
g_{\theta} \in E_{\theta}  \} \ \mbox{ and } \nonumber \\
&&(E_4(A))_{\theta}^{-}=\mathrm{Span}_F\{
a_{\theta} \otimes g_{\theta} | a_{\theta} \in A_{\theta}^{-}, \
g_{\theta} \in E_{\theta} \} \qquad \mbox{ if } \theta \in \{\bar{0}, \bar{1}\};  \\
&&(E_4(A))_{\xi}^{+}=\mathrm{Span}_F\{
a_{\xi} \otimes g_{\xi} | a_{\xi} \in A_{\xi}^{-}, \
g_{\xi} \in E_{\xi} \} \ \mbox{ and } \nonumber \\
&&(E_4(A))_{\xi}^{-}=\mathrm{Span}_F\{
a_{\xi} \otimes g_{\xi} | a_{\xi} \in A_{\xi}^{+}, \
g_{\xi} \in E_{\xi} \} \qquad \mbox{ if } \xi \in \{\bar{2}, \bar{3}\}.  \nonumber
\end{eqnarray}

Let us define some transformations of multilinear $\mathbb{Z}_4$-graded $*$-polynomials.
Denote by $X^{\mathfrak{od}}=Y_{\bar{1}} \cup Z_{\bar{1}} \cup Y_{\bar{3}} \cup Z_{\bar{3}}$ the subset of all variables, odd in respect to the  $\mathbb{Z}_4$-grading, and by $X^{\mathfrak{ev}}=Y_{\bar{0}} \cup Z_{\bar{0}} \cup Y_{\bar{2}} \cup Z_{\bar{2}}$ the subset of all $\mathbb{Z}_4$-even variables.
Fix on $X^{\mathfrak{od}}$ the linear order
$y_{1 \bar{1}} < y_{2 \bar{1}} < \dots < z_{1 \bar{1}} < z_{2 \bar{1}} < \dots < y_{1 \bar{3}} < y_{2 \bar{3}} < \dots < z_{1 \bar{3}} < z_{2 \bar{3}} < \dots$
Assume that $f \in \mathfrak{F}^{(4)}$ is a multilinear graded $*$-polynomial.
Then $f$ is uniquely represented in the form
\begin{equation} \label{oddPol}
f=\sum_{u} \sum_{\sigma \in \mathrm{S}_k} \alpha_{\sigma, u} \ u_1 x_{\sigma(1)} u_2 x_{\sigma(2)} \cdots x_{\sigma(k)} u_{k+1},
\end{equation}
where $x_j \in X^{\mathfrak{od}},$ and $u=u_1 u_2 \cdots u_{k+1}$ is a multilinear monomial over $X^{\mathfrak{ev}},$ possibly empty, \ $k \ge 0.$
Then we assume that
\begin{equation} \label{S}
\mathfrak{s}(f)=\sum_{u} \sum_{\sigma \in \mathrm{S}_k} (-1)^{\sigma} \alpha_{\sigma, u} \  u_1 x_{\sigma(1)} u_2 x_{\sigma(2)} \cdots x_{\sigma(k)} u_{k+1}.
\end{equation}
Consider a collection of variables $(y_{\theta},z_{\theta})=(y_{1 \theta},\dots,y_{n_{\theta} \theta},z_{1 \theta},\dots,z_{m_{\theta} \theta})$ of $\mathbb{Z}_4$-degree $\theta.$
Then for a multilinear graded $*$polynomial $f=f(y_{\bar{0}},z_{\bar{0}},y_{\bar{1}},z_{\bar{1}},y_{\bar{2}},z_{\bar{2}},y_{\bar{3}},z_{\bar{3}})$
\begin{equation} \label{T}
\mathfrak{t}(f)
=f \left|_{\substack{y_{i \bar{2}}:=z_{i \bar{2}}, y_{i \bar{3}}:=z_{i \bar{3}}, \\
 z_{i \bar{2}}:=y_{i \bar{2}}, z_{i \bar{3}}:=y_{i \bar{3}}}} \right.=
f(y_{\bar{0}},z_{\bar{0}},y_{\bar{1}},z_{\bar{1}},z_{\bar{2}},y_{\bar{2}},z_{\bar{3}},y_{\bar{3}})
\end{equation}
is the respective exchange of the variables $y \in Y_{\theta}$ by $z \in Z_{\theta},$ and $z$ by $y$
of $\mathbb{Z}_4$-degrees $\theta=\bar{2}$ and $\bar{3}.$ Observe that
$\mathfrak{t}(y_{1 \theta},\dots,y_{n_{\theta} \theta},z_{1 \theta},\dots,z_{m_{\theta} \theta})=
(z_{1 \theta},\dots,z_{n_{\theta} \theta},y_{1 \theta},\dots,y_{m_{\theta} \theta}).$
It is clear that $\mathfrak{s},$ $\mathfrak{t}$
are linear operators on the space of multilinear $*\mathbb{Z}_4$-polynomials. These operators satisfy
the relations $\mathfrak{s}^2=\mathfrak{t}^2=id,$ \ $\mathfrak{s} \mathfrak{t} = \pm \mathfrak{t} \mathfrak{s},$ where $id$ is the identical transformation, and the sign in the second formula
is defined by the permutation of variables $y_{\bar{3}},z_{\bar{3}}$ induced by applying of $\mathfrak{t}.$
Then we denote
\begin{equation} \label{ST}
\widetilde{f}=\mathfrak{s} \mathfrak{t}(f)
\end{equation}
for a multilinear $*\mathbb{Z}_4$-polynomial $f \in \mathfrak{F}^{(4)}.$ It is clear that
$\widetilde{\widetilde{f}}=\pm f$ for any multilinear $f \in \mathfrak{F}^{(4)}.$ Moreover, we have the next Lemma.

\begin{lemma} \label{EnvId}
A $\mathbb{Z}_4$-graded algebra $A$ with involution satisfies a multilinear
$\mathbb{Z}_4$-graded $*$-identity $f=0$ if and only if $E_4(A)$ satisfies
$\widetilde{f}=0.$
\end{lemma}
\noindent {\bf Proof.}
Assume that $f$ is a multilinear
$\mathbb{Z}_4$-graded $*$-polynomial. Then
\[
f=\sum_{w} \alpha_w \ w\bigl((y_{i_1 \bar{0}}),(z_{i_2 \bar{0}}),(y_{i_3 \bar{1}}),(z_{i_4 \bar{1}}),(y_{i_5 \bar{2}}),(z_{i_6 \bar{2}}),(y_{i_7 \bar{3}}),(z_{i_8 \bar{3}})\bigr), \ \ \alpha_w \in F,\\
\]
where $w=w\bigl((y_{i_1 \bar{0}}),(z_{i_2 \bar{0}}),(y_{i_3 \bar{1}}),(z_{i_4 \bar{1}}),(y_{i_5 \bar{2}}),(z_{i_6 \bar{2}}),(y_{i_7 \bar{3}}),(z_{i_8 \bar{3}})\bigr)$ is a multilinear monomial,
$y_{\theta}=(y_{i \theta}),$ \ $z_{\theta}=(z_{i \theta}),$ \ $\theta \in \mathbb{Z}_4.$
Therefore,
\begin{eqnarray*}
&&\widetilde{f}=\sum_{w} \ (-1)^{\sigma_{\widehat{w}}} \ \alpha_w \ \widehat{w}, \qquad{ where}\\
&&\widehat{w}=\mathfrak{t}(w)=
w\bigl((y_{i_1 \bar{0}}),(z_{i_2 \bar{0}}),(y_{i_3 \bar{1}}),(z_{i_4 \bar{1}}),(z_{i_5 \bar{2}}),(y_{i_6 \bar{2}}),(z_{i_7 \bar{3}}),(y_{i_8 \bar{3}})\bigr)=\\
&&u_{\widehat{w} 1} x_{\sigma_{\widehat{w}}(1)} u_{\widehat{w} 2} x_{\sigma_{\widehat{w}}(2)} \cdots x_{\sigma_{\widehat{w}}(k)} u_{\widehat{w} k+1}.
\end{eqnarray*}
The last formula gives the representation (\ref{oddPol}) of the monomial $\widehat{w};$ here \ $\sigma_{\widehat{w}} \in \mathrm{S}_k,$ \ $x_j \in  X^{\mathfrak{od}},$ and $u_{\widehat{w} j}$ are monomials over $X^{\mathfrak{ev}},$ possibly empty.
Since $\widetilde{f}$ is multilinear then it is enough to consider its evaluations by elements $a \otimes g,$
where $a \in A_{\theta},$ \ $g \in E_{\theta}.$ Taking into account (\ref{Esym})
we need to consider evaluations of the form
\begin{eqnarray} \label{Eval}
  &&y_{i_1 \bar{0}} = b_{i_1 \bar{0}} \otimes h_{i_1 \bar{0}}, \quad  z_{i_2 \bar{0}} = c_{i_2 \bar{0}} \otimes \tilde{h}_{i_2 \bar{0}}, \nonumber\\
  &&y_{i_3 \bar{1}} = b_{i_3 \bar{1}} \otimes g_{i_3 \bar{1}}, \quad  z_{i_4 \bar{1}} = c_{i_4 \bar{1}} \otimes \tilde{g}_{i_4 \bar{1}}, \nonumber\\
  &&y_{i_6 \bar{2}} = c_{i_6 \bar{2}} \otimes h_{i_6 \bar{2}}, \quad  z_{i_5 \bar{2}} = b_{i_5 \bar{2}} \otimes \tilde{h}_{i_5 \bar{2}}, \nonumber\\
  &&y_{i_8 \bar{3}} = c_{i_8 \bar{3}} \otimes g_{i_8 \bar{3}}, \quad  z_{i_7 \bar{3}} = b_{i_7 \bar{3}} \otimes \tilde{g}_{i_7 \bar{3}},
\end{eqnarray}
where  $b_{j \theta} \in A^{+}_{\theta},$ \ $c_{j \theta} \in A^{-}_{\theta},$  and elements \
$h_{j \theta}, g_{j \theta},$ $\tilde{h}_{j \theta}, \tilde{g}_{j \theta} \in E_{\theta}$ \
involve disjoint sets of generators of $E.$
Assume that $(a_1 \otimes g_1,\dots,a_n \otimes g_n)$ is an evaluation of $\widetilde{f}$ of the type (\ref{Eval}) (for corresponding elements $a_i \in A,$ \ $g_i \in E$). Observe that the elements
$h_{i \theta},$ $\tilde{h}_{j \xi} \in E_{\bar{0}} \cup E_{\bar{2}}$ commute with any element of $E,$ and the elements $g_{i \theta},$ $\tilde{g}_{j \xi} \in E_{\bar{1}} \cup E_{\bar{3}}$ anti-commute among themselves.
Then we obtain
\begin{eqnarray*}
&&\widehat{w}(a_1 \otimes g_1,\dots,a_n \otimes g_n) =
w((b_{i_1 \bar{0}} \otimes h_{i_1 \bar{0}}),(c_{i_2 \bar{0}} \otimes \tilde{h}_{i_2 \bar{0}}),(b_{i_3 \bar{1}} \otimes g_{i_3 \bar{1}}),\\ &&(c_{i_4 \bar{1}} \otimes \tilde{g}_{i_4 \bar{1}}),
(b_{i_5 \bar{2}} \otimes \tilde{h}_{i_5 \bar{2}}),(c_{i_6 \bar{2}} \otimes h_{i_6 \bar{2}}),(b_{i_7 \bar{3}} \otimes \tilde{g}_{i_7 \bar{3}}),(c_{i_8 \bar{3}} \otimes g_{i_8 \bar{3}}))=\\
&&w((b_{i_1 \bar{0}}),(c_{i_2 \bar{0}}),(b_{i_3 \bar{1}}),(c_{i_4 \bar{1}}),(b_{i_5 \bar{2}}),
(c_{i_6 \bar{2}}),(b_{i_7 \bar{3}}),(c_{i_8 \bar{3}})) \ \otimes \ \widehat{w}(g_1,\dots,g_n)=\\
&&w(\tilde{a}_1,\dots,\tilde{a}_n) \ \otimes \  u_{\widehat{w} 1}(h,\tilde{h}) \
g'_{\sigma_{\widehat{w}}(1)} \  u_{\widehat{w} 2}(h,\tilde{h}) \ g'_{\sigma_{\widehat{w}}(2)} \ \cdots \ g'_{\sigma_{\widehat{w}}(k)} \
u_{\widehat{w} k+1}(h,\tilde{h})=\\
&&(-1)^{\sigma_{\widehat{w}}} \
w(\tilde{a}_1,\dots,\tilde{a}_n) \ \otimes \ g_1 \cdots g_n.
\end{eqnarray*}
Where $(\tilde{a}_1,\dots,\tilde{a}_n)=((b_{i_1 \bar{0}}),(c_{i_2 \bar{0}}),(b_{i_3 \bar{1}}),(c_{i_4 \bar{1}}),(b_{i_5 \bar{2}}), (c_{i_6 \bar{2}}),(b_{i_7 \bar{3}}),(c_{i_8 \bar{3}}))$ are arbitrary $\widehat{G}$-homogeneous elements of $A,$ \   $u_{\widehat{w} j}(h,\tilde{h})$ are monomials $u_{\widehat{w} j}$ evaluated by elements $h_{i \theta},$
$\tilde{h}_{j \xi} \in E_{\bar{0}} \cup E_{\bar{2}},$ and the $k$-tuple $(g'_1,\dots,g'_k)=((g_{i_3 \bar{1}}),(\tilde{g}_{i_4 \bar{1}}),(\tilde{g}_{i_7 \bar{3}}),(g_{i_8 \bar{3}})).$
Therefore,
\begin{eqnarray*}
&&\widetilde{f}(a_1 \otimes g_1,\dots,a_n \otimes g_n) =\sum_{w} \ (-1)^{\sigma_{\widehat{w}}} \  \alpha_w \
\widehat{w}(a_1 \otimes g_1,\dots,a_n \otimes g_n) = \\
&& \sum_{w} \ (-1)^{\sigma_{\widehat{w}}} (-1)^{\sigma_{\widehat{w}}} \  \alpha_w \
w(\tilde{a}_1,\dots,\tilde{a}_n) \ \otimes \ g_1 \cdots g_n
=f(\tilde{a}_1,\dots,\tilde{a}_n) \otimes g_1 \cdots g_n.
\end{eqnarray*}
Thus $\widetilde{f}(a_1 \otimes g_1,\dots,a_n \otimes g_n)=0$ for any evaluation (\ref{Eval}) if and only if $f(\tilde{a}_1,\dots,\tilde{a}_n)=0$ for all appropriate $\tilde{a}_i \in A^{\delta_i}_{\theta_i},$ \
$\delta_i \in \{ +, - \},$ \ $\theta_i \in \mathbb{Z}_4.$
\hfill $\Box$

\begin{definition}
Given a $gi$T-ideal $\Gamma \subseteq \mathfrak{F}^{(4)}$ denote by $\widetilde{\Gamma}$
the $gi$T-ideal generated by the set $S=\{ \widetilde{f} | f \in \Gamma \cap (\cup_{n \ge 1} P_n) \ \}$
of $\mathfrak{s} \mathfrak{t}$-images of all multilinear polynomials from $\Gamma.$
\end{definition}

Lemma \ref{EnvId} along with properties of the operators $\mathfrak{s},$ $\mathfrak{t}$
immediately implies the following.

\begin{lemma} \label{EnvId1}
Given a $gi$T-ideal $\Gamma \subseteq \mathfrak{F}^{(4)}$ we have that $\Gamma=\mathrm{Id}^{gi}(A)$ for a $\mathbb{Z}_4$-graded $*$-algebra $A$ iff \ $\widetilde{\Gamma}=\mathrm{Id}^{gi}(E_4(A)).$ Besides that, $\widetilde{\widetilde{\Gamma}}=\Gamma.$
\end{lemma}

Hence, we have that $A \sim_{gi} B$ for $\mathbb{Z}_4$-graded $*$-algebras $A,$
$B$ if and only if $E_4(A) \sim_{gi} E_4(B).$ And $E_4(E_4(A)) \sim_{gi} A$
for any $\mathbb{Z}_4$-graded algebra $A$ with involution. The last property is
also a simple consequence of the facts that
$E_4(E_4(A))=\bigoplus_{\theta \in \mathbb{Z}_4} A_{\theta} \otimes_F E_{\theta} \otimes_F E_{\theta},$
and the algebra $E_4(E)=\bigoplus_{\theta \in \mathbb{Z}_4} E_{\theta} \otimes_F E_{\theta}$ is commutative and non-nilpotent.

\begin{remark} \label{Z4Var}
Since $E_4(A)$ is a subalgebra of $A \otimes_F E$ then by Regev's theorem \cite{Reg1} we have that
$E_4(A)$ is a PI-algebra if and only if $A$ is a PI-algebra.
Particularly, consider a $*$-variety $\mathfrak{V}.$ Assume that $\mathfrak{V}$ is defined by a $*$T-ideal $\Gamma \subseteq F\langle Y, Z \rangle,$ and $\Gamma=\mathrm{Id}^{*}(A)$ for an algebra with involution $A.$
Denote by $\widetilde{\mathfrak{V}}^{\mathbb{Z}_4}$ the class of all associative $\mathbb{Z}_4$-graded $F$-algebras $B$ with involution such that $E_4(B) \in \mathfrak{V}.$ It is clear from Lemma \ref{EnvId1}
that $\widetilde{\mathfrak{V}}^{\mathbb{Z}_4}$ is a $\mathbb{Z}_4$-graded $*$-variety defined
by the $gi$T-ideal $\Gamma_1$ of $\mathbb{Z}_4$-graded $*$-identities of the $\mathbb{Z}_4$-graded
algebra with involution $A \otimes_F E=\bigoplus_{\theta \in \mathbb{Z}_4} A \otimes_F E_{\theta}.$
The $gi$T-ideal $\Gamma_1=\widetilde{\Gamma}_2,$ where $\Gamma_2$
is the $gi$T-ideal generated by $\Gamma,$ i.e.
\begin{equation} \label{giT*}
\Gamma_2=\Gamma^{\mathbb{Z}_4}=
giT[S_{\Gamma}] \quad  \mbox{for } \ \
S_{\Gamma}=\{ \  f |_{y_{i}:=\sum_{\theta \in \mathbb{Z}_4} y_{i \theta}, \
z_{i}:=\sum_{\theta \in \mathbb{Z}_4} z_{i \theta}, \  \forall i}  \  | \ f \in \Gamma \ \}.
\end{equation}
\end{remark}

\section{Alternating and symmetrizing polynomials.}

Let $f=f(s_1,\dots,s_k,x_1,\dots,x_n) \in F\langle Y, Z \rangle$
be a multilinear polynomial. Assume that $S=\{ s_1,
\dots, s_k\} \subseteq Y$ or $S \subseteq Z$. We say
that $f$ is alternating in $S,$ if
$f(s_{\sigma(1)},\dots,s_{\sigma(k) },x_1,\dots,x_n)=(-1)^{\sigma}
f(s_1,\dots,s_k,x_1,\dots,x_n)$ holds for any permutation $\sigma
\in \mathrm{S}_k.$

For any multilinear polynomial with involution
$g(s_1,\dots,s_k,x_1,\dots,x_n)$ we construct a multilinear polynomial $f$ alternating in
$S=\{ s_1, \dots, s_k\}$ by setting \[f(s_1,\dots,s_k,x_1,\dots,x_n)=
\mathcal{A}_{S}(g)=\sum_{\sigma \in \mathrm{S}_k} (-1)^{\sigma}
g(s_{\sigma(1)},\dots,s_{\sigma(k)},x_1,\dots,x_n).\] The
corresponding mapping $\mathcal{A}_{S}$ is a linear transformation of multilinear $*$-po\-ly\-no\-mials.
We call it the alternator. Any $*$-polynomial $f$ alternating in $S$
can be decomposed as $f=\sum_{i=1}^{m} \alpha_i
\mathcal{A}_{S}(u_i),$ where the $u_i$'s are monomials, \
$\alpha_i \in F.$

We say that a multilinear $*$-polynomial $f(s_1,\dots,s_k,x_1,\dots,x_n)$ is symmetrizing in the set
$S=\{ s_1,\dots,s_k \}$ ($S \subseteq Y$ or $S \subseteq Z$), if
$f(s_{\sigma(1)},\dots,s_{\sigma(k)},x_1,\dots,x_n)=
f(s_1,\dots,s_k,x_1,\dots,x_n)$ for any $\sigma \in
\mathrm{S}_k.$

For any multilinear $*$-polynomial
$g(s_1,\dots,s_k,x_1,\dots,x_n)$ the multilinear $*$-po\-ly\-no\-mial
$$f(s_1,\dots,s_k,x_1,\dots,x_n)=\mathcal{E}_{S}(g)=\sum_{\sigma
\in \mathrm{S}_k} g(s_{\sigma(1)},\dots,s_{\sigma(k)},x_1,\dots,x_n). $$
is symmetrizing in $S.$ $\mathcal{E}_{S}$ is also a linear
transformation of multilinear $*$-polynomials. It is called the symmetrizator.
Any multilinear $*$-polynomial $f$ symmetrizing in $S$ can be written as
$f=\sum_{i=1}^{m} \alpha_i \mathcal{E}_{S}(u_i),$ where the
$u_i$'s are monomials of $f,$ and $\alpha_i \in F.$
Properties of alternating and symmetrizing polynomials with
involution are similar to that of ordinary polynomials (see,
e.g., \cite{Dren}, \cite{GZbook}, \cite{Kem1}).

Particularly, a multilinear $*$-polynomial $f(s_1,\dots,s_k,x_1,\dots,x_n)$
is symmetrizing in variables $S=\{ s_1,\dots,s_k \}$ iff $f(s_1,\dots,s_k,x_1,\dots,x_n)$
is the full linearization in the variable $s$ of the non-zero polynomial $\hat{f}=\frac{1}{k!} f(s,\dots,s,x_1,\dots,x_n).$
Moreover, $\hat{f}=\sum_{(v)} \alpha_{(v)} v_{0} s v_{1} s \cdots s v_{k},$
whenever $f=\sum_{(v)} \sum_{\sigma
\in \mathrm{S}_k} \alpha_{(v)} v_{0} s_{\sigma(1)} v_{1} s_{\sigma(1)} \cdots s_{\sigma(k)} v_{k},$
where monomials $v_{0}, v_{1}, \dots, v_{k}$ (possibly empty) do not depend on $S,$ \ 
$\alpha_{(v)} \in F.$

Similarly we can consider graded $*$-polynomials alternating or symmetrizing in
a set of variables $S \subseteq Y_{\theta}$ or $S \subseteq Z_{\theta}$ for any
fixed $\theta \in G$ (see \cite{Svi3}).

\begin{lemma} \label{SProp}
Consider disjoint collections of variables
$\bar{y}=\{ y_1,\dots,y_n \} \subseteq Y_{\bar{0}} \cup Y_{\bar{1}},$ \
$\bar{z}=\{ z_1,\dots,z_m \} \subseteq Z_{\bar{0}} \cup Z_{\bar{1}},$ and
$\bar{t}=\{ t_{1 1},\dots,t_{1 \hat{n}_1},\dots,t_{\hat{k} 1},\dots,t_{\hat{k} \hat{n}_{\hat{k}}} \},$
where $\{ t_{i 1},\dots,t_{i \hat{n}_i} \} \subseteq Y_{\bar{1}},$ or
$\{ t_{i 1},\dots,t_{i \hat{n}_i} \}  \subseteq Z_{\bar{1}}$ for any $i=1,\dots,\hat{k}.$
Let $f(\bar{y},\bar{z},\bar{t}) \in \mathfrak{F}^{(4)}$ be a multilinear graded
$*$-polynomial, which is alternating in any collection $\{ t_{i 1},\dots,t_{i \hat{n}_i} \},$ \  $i=1,\dots,\hat{k}.$ Then the polynomial $\widetilde{f}$ depends on the same variables as $f,$
and $\widetilde{f}$ is symmetrizing in $\{ t_{i 1},\dots,t_{i \hat{n}_i} \}$ for any $i=1,\dots,\hat{k}.$ \
\end{lemma}
\noindent {\bf Proof.}
It is clear that $\widetilde{f}=\mathfrak{s}(f).$ Also the polynomial $f$ can be decomposed as
\[f=\sum_{\substack{(v), \\ \tau \in \mathrm{S}_r}} \sum_{\substack{\sigma_i
\in \mathrm{S}_{\hat{n}_i},\\ {\scriptscriptstyle 1 \le i \le \hat{k}} }}
\alpha_{(v), \tau} \ (-1)^{\sigma_1} \cdots (-1)^{\sigma_{\hat{k}}} \
(\sigma_1 \cdots \sigma_{\hat{k}} \  \tau) \  v_{0}  x_1  v_{1}  x_2 \cdots  x_r  v_{k},\]
where $x_j \in Y_{\bar{1}} \cup Z_{\bar{1}},$ $v_j$ are monomials (possibly empty) over
$Y_{\bar{0}} \cup Z_{\bar{0}},$ the permutation $\tau$ acts on the variables $x_j,$ and
the permutations $\sigma_i$ acts on disjoint subsets of the set $\{ x_1, \dots, x_r \}$ corresponding to
the sets of variables $\{ t_{i 1},\dots,t_{i \hat{n}_i} \},$ \ $\alpha_{(v), \tau} \in F.$ Then
\begin{eqnarray*}
&&\widetilde{f}=\sum_{\substack{(v), \\ \tau \in \mathrm{S}_r}}
\sum_{\substack{\sigma_i
\in \mathrm{S}_{\hat{n}_i},\\ {\scriptscriptstyle 1 \le i \le \hat{k}} }}
 \alpha_{(v), \tau}
(-1)^{\sigma_1 \cdots \sigma_{\hat{k}}} \
(-1)^{\sigma_1 \cdots \sigma_{\hat{k}} \tau} \  (\sigma_1 \cdots \sigma_{\hat{k}} \  \tau) \  v_{0} x_1 v_{1}  x_2  \cdots  x_r  v_{k} =\\
&&\sum_{\substack{(v), \\ \tau \in \mathrm{S}_r}} (-1)^{\tau} \alpha_{(v), \tau}
\sum_{\substack{\sigma_i
\in \mathrm{S}_{\hat{n}_i},\\ {\scriptscriptstyle 1 \le i \le \hat{k}} }} \
(\sigma_1 \cdots \sigma_{\hat{k}}) \  v_{0} x_{\tau(1)}  v_{1}  x_{\tau(2)}  \cdots  x_{\tau(r)}  v_{k}.
\end{eqnarray*}
Thus $\widetilde{f}$ is symmetrizing in any $\{ t_{i 1},\dots,t_{i \hat{n}_i} \}.$
\hfill $\Box$

Given a $*$T-ideal $\Gamma \unlhd F\langle Y, Z \rangle$
the vector space $\Gamma_{n,m}=\Gamma \cap P_{n,m}$
of multilinear $*$-polynomials $f(y_1,\dots,y_{n},z_1,\dots,z_{m}) \in \Gamma$
has the structure of $(FS_{n} \otimes FS_{m})$-module defined by
$(\sigma \otimes \tau) f(y_1,\dots,y_{n},z_1,\dots,z_{m})=
f(y_{\sigma(1)},\dots,y_{\sigma(n)},z_{\tau(1)},\dots,z_{\tau(m)})$ for any
$(\sigma,\tau) \in S_{n} \times S_{m}.$ The character of the quotient module
$P_{n,m}(\Gamma)=P_{n,m}/\Gamma_{n,m} \subseteq F\langle Y, Z \rangle/\Gamma$
can be decomposed as $\chi_{n,m}(\Gamma)=\sum_{\substack{\lambda \vdash n\\ \mu \vdash m}} \mathfrak{m}_{\lambda, \mu} \  (\chi_{\lambda} \otimes \chi_{\mu}),$
where $\chi_{\lambda} \otimes \chi_{\mu}$ is the irreducible $S_{n} \times S_{m}$-character
associated to the pair $(\lambda, \mu)$ of partitions $\lambda \vdash n,$ $\mu \vdash m,$ \
$\mathfrak{m}_{\lambda, \mu} \in \mathbb{Z}$ is a multiplicity (see for instance \cite{GM}, \cite{GM1}, \cite{GZbook}, \cite{JKer}). An irreducible submodule of $P_{n,m}(\Gamma)$
corresponding to the pair $(\lambda, \mu)$ is generated by a non-zero polynomial
$f_{\lambda, \mu}=(e_{T_{\lambda}} \otimes e_{T_{\mu}}) f,$ where $f \in P_{n,m},$
and $e_{T_{\lambda}} \in FS_{n},$ \  $e_{T_{\mu}} \in FS_{m}$ are the essential idempotents
corresponding to the Young tableaux $T_{\lambda},$ and $T_{\mu}$ respectively
(see Definition 2.2.12 \cite{GZbook}). We say that a multilinear $*$-polynomial $f$
corresponds to the pair of partitions $(\lambda, \mu)$ if
$(FS_{n} \otimes FS_{m}) \  f = (FS_{n} \otimes FS_{m}) \  f_{\lambda, \mu}.$
Particularly, the next observation holds.

\begin{remark} \label{Repr}
Given a multilinear $*$-polynomial $f \in P_{n,m}$ there exist a finite set of pairs $(\lambda_j, \mu_j)$ (not necessary different) of partitions $\lambda_j \vdash n,$ $\mu_j \vdash m$ ($j=1,\dots,k$) and multilinear
$*$-polynomials $g_{\lambda_j, \mu_j} \in P_{n,m}$ such that $g_{\lambda_j, \mu_j}$ corresponds to $(\lambda_j, \mu_j),$ and the $*$T-ideal generated by $f$ can be decomposed as $*T[f]=\sum_{j=1}^{k} *T[g_{\lambda_j, \mu_j}].$
\end{remark}

Moreover, by Theorem 5.9 \cite{GR}
\begin{equation} \label{hook}
\chi_{n,m}(\Gamma)=\sum_{(\lambda,\mu) \in H_{\Gamma}}
\mathfrak{m}_{\lambda, \mu} \  (\chi_{\lambda} \otimes \chi_{\mu}),
\end{equation}
where $H_{\Gamma}=(H(k_1,l_1),H(k_2,l_2))$ is a double hook corresponding to $\Gamma.$
The hook $H(k,l)$ is the set of all partitions $\lambda=(\lambda_1,\dots,\lambda_s)$ satisfying
the condition $\lambda_{k+1} \le l.$ Applying arguments of Lemma 2.5.6 \cite{GZbook} we always can assume
that for any $(\lambda,\mu) \in H_{\Gamma}$ the set of variables of a polynomial $f_{\lambda, \mu}$
can be decomposed into disjoint unions $\{ y_1, \dots, y_n \}=(\bigcup_{i=1}^{k'_1} Y'_i) \bigcup
(\bigcup_{i=1}^{l'_1} T'_i),$ \ $\{ z_1, \dots, z_m \}=(\bigcup_{i=1}^{k'_2} Z'_i) \bigcup
(\bigcup_{i=1}^{l'_2} S'_i),$ where $k'_r \le k_r,$ \ $l'_r \le l_r$ ($r=1,2$), and
$f_{\lambda, \mu}$ is symmetrizing in any $Y'_i \subseteq Y,$ $Z'_j \subseteq Z$ ($1 \le i \le k'_1,$
$1 \le j \le k'_2$), and alternating in any $T'_i \subseteq Y,$ $S'_j \subseteq Z$
($1 \le i \le l'_1,$ $1 \le j \le l'_2$).
Notice that $\mathfrak{m}_{\lambda, \mu}=0$ in \eqref{hook}
means that $f_{\lambda, \mu}=(e_{T_{\lambda}} \otimes e_{T_{\mu}}) f \in \Gamma$
for any Young tableaux $T_{\lambda},$ $T_{\mu}$ and for any $*$-polynomial $f \in P_{n,m},$
(see, e.g., Theorem 2.4.5 \cite{GZbook}).

\section{Classification theorems.}

\begin{theorem} \label{IPI1}
Let $F$ be a field of characteristic zero.
Any proper $*$T-ideal of the free associative $F$-algebra with involution
is the ideal of identities with involution of the Grassmann
$\mathbb{Z}_4$-envelope of some finitely generated associative $\mathbb{Z}_4$-graded
PI-algebra with graded involution.
\end{theorem}
\noindent {\bf Proof.}
Let $\Gamma$ be a proper $*$T-ideal of $F\langle Y, Z \rangle,$ and $\mathfrak{V}_{\Gamma}$ the $*$-variety
defined by $\Gamma.$ Consider the $\mathbb{Z}_4$-graded $*$-variety $\widetilde{\mathfrak{V}}_{\Gamma}^{\mathbb{Z}_4}$ of all associative $\mathbb{Z}_4$-graded $*$-algebras $B$ such that $E_4(B) \in \mathfrak{V}_{\Gamma}.$
Assume that $H_{\Gamma}=(H(k_1,l_1),H(k_2,l_2))$ is the double hook corresponding to $\Gamma$ \cite{GR}. Take $\nu=\max\{k_1,l_1,k_2,l_2\},$ and the relatively free algebra $\mathcal{R}$ of the rank $\nu$
of the $\mathbb{Z}_4$-graded $*$-variety $\widetilde{\mathfrak{V}}_{\Gamma}^{\mathbb{Z}_4}.$
Then as in Remark \ref{Z4Var} we have $\mathcal{R}=
F\langle Y_{\nu}^{\mathbb{Z}_4}, Z_{\nu}^{\mathbb{Z}_4} \rangle /(\widetilde{\Gamma}_2 \cap F\langle Y_{\nu}^{\mathbb{Z}_4}, Z_{\nu}^{\mathbb{Z}_4} \rangle ),$ where $\Gamma_2=\Gamma^{\mathbb{Z}_4}=
giT[S_{\Gamma}]$ is defined by (\ref{giT*}).
By Remark \ref{Z4Var} and Amitsur's theorem \cite{A3}, \cite{A4} $\mathcal{R}$ is a PI-algebra.
Let us prove that $E_4(\mathcal{R})$ generates $\mathfrak{V}_{\Gamma}.$

It is clear that
$\mathrm{Id}^{*}(E_4(\mathcal{R})) \supseteq \Gamma.$ Take a multilinear polynomial with involution
$f(y_1,\dots,y_n,z_1,\dots,z_m) \in \mathrm{Id}^{*}(E_4(\mathcal{R})) \cap P_{n,m}.$ By Remark \ref{Repr} we can assume that $f$ corresponds to a pair of partitions $(\lambda,\mu),$ where $\lambda \vdash n,$ and $\mu \vdash m.$ If $(\lambda,\mu) \notin H_{\Gamma}$ then $f \in \Gamma$ by Theorem 5.9 \cite{GR}. Suppose that
$(\lambda,\mu) \in H_{\Gamma}.$ Then similarly to Lemma 2.5.6 \cite{GZbook} we can assume that the set $\{ y_1,\dots,y_n \} \in Y$ of the variables of $f$ is divided on at most $\nu$ sets of symmetrized variables
$\{ y_{i_1 1},\dots,y_{i_1 n_{i_1}} \}$ ($i_1 =1,\dots,\nu$), and  at most
$\nu$ sets of alternated variables $\{ t_{i_2 1},\dots,t_{i_2 \hat{n}_{i_2}} \}$ ($i_2 =1,\dots,\nu$).
Similarly, the set $\{ z_1,\dots,z_m \}$ consists of at most $\nu$ sets of symmetrized variables
$\{ z_{j_1 1},\dots,z_{j_1 m_{j_1}} \}$ ($j_1 =1,\dots,\nu$), and  at most
$\nu$ sets of alternated variables $\{ s_{j_2 1},\dots,s_{j_2 \hat{m}_{j_2}} \}$ ($j_2 =1,\dots,\nu$).
Thus, $f=f(\vec{y},\vec{t},\vec{z},\vec{s}),$
where
\begin{eqnarray} \label{Var}
&&\vec{y}=(y_{1 1},\dots,y_{1 n_1},\dots,y_{\nu 1},\dots,y_{\nu n_{\nu}}) \subseteq Y, \nonumber\\
&&\vec{t}=(t_{1 1},\dots,t_{1 \hat{n}_1},\dots,t_{\nu 1},\dots,t_{\nu \hat{n}_{\nu}}) \subseteq Y, \nonumber\\
&&\vec{z}=(z_{1 1},\dots,z_{1 m_1},\dots,z_{\nu 1},\dots,z_{\nu m_{\nu}}) \subseteq Z, \nonumber\\
&&\vec{s}=(s_{1 1},\dots,s_{1 \hat{m}_1},\dots,z_{\nu 1},\dots,z_{\nu \hat{m}_{\nu}}) \subseteq Z
\end{eqnarray}
are disjoint collections of variables, and $f$ is symmetrizing in any $\{ y_{i 1},\dots,y_{i n_{i}} \},$ and $\{ z_{i 1},\dots,z_{i m_{i}} \},$ and alternating in any $\{ t_{i 1},\dots,t_{i \hat{n}_{i}} \},$ and $\{ s_{i 1},\dots,s_{i \hat{m}_{i}} \}$ ($i=1,\dots,\nu$).

Since $f \in \mathrm{Id}^{*}(E_4(\mathcal{R}))$ then $f$ is equal to zero in $E_4(\mathcal{R})$ for
\begin{eqnarray} \label{VarEv}
&&y_{i j_1} = \bar{y}_{i \bar{0}} \otimes h_{n \cdot i +j_1 \,  \bar{0}},
\qquad
t_{i j_2} = \bar{y}_{i \bar{1}} \otimes g_{n \cdot i +j_2 \,  \bar{1}}, \nonumber\\
&&z_{i j_3} = \bar{z}_{i \bar{0}} \otimes \tilde{h}_{m \cdot i +j_3 \,  \bar{0}}, \qquad
s_{i j_4} = \bar{z}_{i \bar{1}} \otimes \tilde{g}_{m \cdot i +j_4 \,  \bar{1}},\\
&& i=1, \dots, \nu, \quad 1 \le j_1 \le n_{i}, \ 1 \le j_2 \le \hat{n}_{i}, \  1 \le j_3 \le m_{i}, \
1 \le j_4 \le \hat{m}_{i}, \nonumber
\end{eqnarray}
where $\bar{y}_{i \theta}=y_{i \theta} + \mathcal{I},$ \ $\bar{z}_{i \theta}=z_{i \theta} + \mathcal{I},$ \ $y_{i \theta} \in Y_{\theta},$ and $z_{i \theta} \in Z_{\theta}$ are graded variables from
$Y_{\nu}^{\mathbb{Z}_4} \cup Z_{\nu}^{\mathbb{Z}_4}$ of $\mathbb{Z}_4$-degree $\theta \in \{ \bar{0}, \bar{1}\},$ \  $\mathcal{I}=\widetilde{\Gamma}_2 \cap F\langle Y_{\nu}^{\mathbb{Z}_4}, Z_{\nu}^{\mathbb{Z}_4} \rangle,$ \ $h_{l \bar{0}},$ $\tilde{h}_{l \bar{0}} \in E_{\bar{0}},$ \
$g_{l \bar{1}},$ $\tilde{g}_{l \bar{1}} \in E_{\bar{1}}$ are elements of the Grassmann algebra depending on
disjoint sets of generators. Let us denote $a^{(k)}= \underbrace{a,\dots,a}_{k}$ for any element $a.$
Therefore, we obtain in the algebra $E_4(\mathcal{R})$ the equalities
\begin{eqnarray} \label{subst}
&&f|_{(\ref{VarEv})}=\bar{f}_3 \otimes g =0, \qquad \mbox{ where }  \\
&&\bar{f}_3=f_2(\bar{y}_{1 \bar{0}}^{(n_1)},\dots,\bar{y}_{\nu \bar{0}}^{(n_{\nu})},\bar{y}_{1 \bar{1}}^{(\hat{n}_1)},\dots,\bar{y}_{\nu \bar{1}}^{(\hat{n}_{\nu})},\bar{z}_{1 \bar{0}}^{(m_1)},\dots,\bar{z}_{\nu \bar{0}}^{(m_{\nu})},\bar{z}_{1 \bar{1}}^{(\hat{m}_1)},\dots,\bar{z}_{\nu \bar{1}}^{(\hat{m}_{\nu})}). \nonumber
\end{eqnarray}
Here $f_2=\widetilde{f}_1.$ Where the graded multilinear polynomial
$f_1=f(\vec{y}_{\bar{0}},\vec{y}_{\bar{1}},\vec{z}_{\bar{0}},\vec{z}_{\bar{1}}),$  with
\begin{eqnarray} \label{Var2}
&&\vec{y}_{\bar{0}}=(y_{(1,1) \bar{0}},\dots,y_{(1,n_1) \bar{0}},\dots,y_{(\nu,1) \bar{0}},\dots,y_{(\nu, n_{\nu}) \bar{0}}) \subseteq Y_{\bar{0}}, \nonumber\\
&&\vec{y}_{\bar{1}}=(y_{(1,1) \bar{1}},\dots,y_{(1,\hat{n}_1) \bar{1}},\dots,
y_{(\nu,1) \bar{1}},\dots,y_{(\nu,\hat{n}_{\nu}) \bar{1}}) \subseteq Y_{\bar{1}}, \nonumber\\
&&\vec{z}_{\bar{0}}=(z_{(1,1) \bar{0}},\dots,z_{(1,m_1) \bar{0}},\dots,z_{(\nu,1) \bar{0}},\dots,z_{(\nu,m_{\nu}) \bar{0}}) \subseteq Z_{\bar{0}}, \nonumber\\
&&\vec{z}_{\bar{1}}=(z_{(1,1) \bar{1}},\dots,z_{(1,\hat{m}_1) \bar{1}},\dots,z_{(\nu,1) \bar{1}},\dots,z_{(\nu,\hat{m}_{\nu}) \bar{1}}) \subseteq Z_{\bar{1}},
\end{eqnarray}
is the result of the evaluation of the variables $y_{i j},$ $z_{i j}$ of the polynomial $f$
by the corresponding graded variables of the degree $\bar{0},$ and of the variables $t_{i j},$ $s_{i j}$
by the graded variables $y_{(i,j) \bar{1}},$ $z_{(i,j) \bar{1}}$  of the degree $\bar{1}$ respectively.
The element $g$ in \eqref{subst} is the product of all elements $h_{l \bar{0}},$ $\tilde{h}_{l \bar{0}},$ \
$g_{l \bar{1}},$ $\tilde{g}_{l \bar{1}}$ of the Grassmann algebra from (\ref{VarEv}).

Observe that by Lemma \ref{SProp} the polynomial $f_2=\widetilde{f}_1$ is symmetrizing in any set of variables $y_{(i,1) \theta},\dots,y_{(i,n'_i) \theta},$ and $z_{(i,1) \theta},\dots,z_{(i,m'_i) \theta},$ \
for all $i=1,\dots,\nu,$ \ $\theta \in \{ \bar{0}, \bar{1} \}$ (if $\theta=\bar{0}$ then $n'_i=n_i,$
$m'_i=m_i,$ otherwise $n'_i=\hat{n}_i,$ \  $m'_i=\hat{m}_i$).
The equality (\ref{subst}) means that the graded $*$-polynomial
\[f_3=f_2(y_{1 \bar{0}}^{(n_1)},\dots,y_{\nu \bar{0}}^{(n_{\nu})},y_{1 \bar{1}}^{(\hat{n}_1)},\dots,
y_{\nu \bar{1}}^{(\hat{n}_{\nu})},z_{1 \bar{0}}^{(m_1)},\dots,z_{\nu \bar{0}}^{(m_{\nu})},z_{1 \bar{1}}^{(\hat{m}_1)},\dots,z_{\nu \bar{1}}^{(\hat{m}_{\nu})})\]
belongs to $\widetilde{\Gamma}_2 \cap F\langle Y_{\nu}^{\mathbb{Z}_4}, Z_{\nu}^{\mathbb{Z}_4} \rangle.$
Thus, $f_3 \in \widetilde{\Gamma}_2.$ The polynomial
$f_2=\widetilde{f}_1(\vec{y}_{\bar{0}},\vec{y}_{\bar{1}},\vec{z}_{\bar{0}},\vec{z}_{\bar{1}})$
is the full linearization of $\frac{1}{\alpha} \cdot f_3,$ where $\alpha \in F$ is some nonzero coefficient
which appears as the result of identifying of symmetrized variables. The variables of $f_2$ as in (\ref{Var2}). Hence $f_2 \in \widetilde{\Gamma}_2.$

Take the relatively free $*$-algebra $\mathcal{L}=F\langle Y, Z \rangle /\Gamma$ of the $*$-variety
$\mathfrak{V}_{\Gamma},$ and consider the $\mathbb{Z}_4$-graded $*$-algebra
$\mathcal{L} \otimes E=\bigoplus_{\theta \in \mathbb{Z}_4} \mathcal{L} \otimes E_{\theta}.$
By Remark \ref{Z4Var} $\mathcal{L} \otimes E$ satisfies the graded $*$-identity $f_2(\vec{y}_{\bar{0}},\vec{y}_{\bar{1}},\vec{z}_{\bar{0}},\vec{z}_{\bar{1}})=0.$ Particularly, the evaluation
\begin{eqnarray} \label{VarEv2}
&&y_{(i,j_1) \bar{0}} = \bar{y}_{i j_1} \otimes h_{n \cdot i +j_1 \, \bar{0}},
\qquad
y_{(i,j_2) \bar{1}} = \bar{t}_{i j_2} \otimes g_{n \cdot i +j_2  \, \bar{1}}, \nonumber\\
&&z_{(i,j_3) \bar{0}} = \bar{z}_{i j_3} \otimes \tilde{h}_{m \cdot i +j_3 \, \bar{0}},
\qquad
z_{(i,j_4) \bar{1}} = \bar{s}_{i j_4} \otimes \tilde{g}_{m \cdot i +j_4 \, \bar{1}},  \\
&& i=1, \dots, \nu, \quad 1 \le j_1 \le n_{i}, \ 1 \le j_2 \le \hat{n}_{i}, \  1 \le j_3 \le m_{i}, \
1 \le j_4 \le \hat{m}_{i} \nonumber
\end{eqnarray}
gives the result
$f_2 |_{(\ref{VarEv2})}=
\widetilde{\widetilde{f}}_1(\vec{\bar{y}},\vec{\bar{t}},\vec{\bar{z}},\vec{\bar{s}}) \otimes g' = f(\vec{\bar{y}},\vec{\bar{t}},\vec{\bar{z}},\vec{\bar{s}}) \otimes g' =0.$
Here $\vec{\bar{y}},\vec{\bar{t}},\vec{\bar{z}},\vec{\bar{s}}$ is the sequence formed as in (\ref{Var}) by the elements $\bar{y}_{i j_1}=y_{i j_1}+\Gamma,$ \ $\bar{t}_{i j_2}=t_{i j_2}+\Gamma,$ \
$\bar{z}_{i j_3}=z_{i j_3}+\Gamma,$ \ $\bar{s}_{i j_4}=s_{i j_4}+\Gamma$ ($i=1, \dots, \nu,$ \quad $1 \le j_1 \le n_{i},$ \ $1 \le j_2 \le \hat{n}_{i},$ \  $1 \le j_3 \le m_{i},$ \ $1 \le j_4 \le \hat{m}_{i}$), where
the variables $y_{i j_1},$ $t_{i j_2},$ $z_{i j_3},$ $s_{i j_4}$ are the same as in (\ref{Var}). The element $g'$ is the product of all elements of the Grassmann algebra from (\ref{VarEv2}) depending on disjoint sets of generators. Therefore, $f(\vec{\bar{y}},\vec{\bar{t}},\vec{\bar{z}},\vec{\bar{s}})=0$ in $\mathcal{L},$
and $f \in \Gamma.$ Hence $\mathrm{Id}^{*}(E_4(\mathcal{R})) = \Gamma.$
\hfill $\Box$

We can reinforce the result similarly to the classical case of Kemer's theorems
for PI-algebras \cite{Kem1} using Theorem 6.2 \cite{Svi3}.

\begin{theorem} \label{IPI2}
Let $F$ be a field of characteristic zero.
Any proper $*$T-ideal of the free associative $F$-algebra with involution
is the ideal of identities with involution of the Grassmann
$\mathbb{Z}_4$-envelope of some associative $\mathbb{Z}_4$-graded
algebra with graded involution, finite dimensional over $F$.
\end{theorem}
\noindent {\bf Proof.} If $\Gamma$ is a proper $*$T-ideal of $F\langle Y, Z \rangle$ then
by Theorem \ref{IPI1} we have $\Gamma=\mathrm{Id}^{*}(E_4(B))$ for some associative finitely
generated $\mathbb{Z}_4$-graded PI-algebra $B$ with graded involution. Theorem 6.2 \cite{Svi3}
states that there exists a finite dimensional over $F$ $\mathbb{Z}_4$-graded algebra $C$ with graded involution which has the same graded $*$-identities as $B.$ Hence $E_4(B) \sim_{gi} E_4(C).$
Particularly, $\mathrm{Id}^{*}(E_4(B))=\mathrm{Id}^{*}(E_4(C))=\Gamma.$
\hfill $\Box$

For a finitely generated associative PI-algebra with involution we also have
the next theorem.

\begin{theorem} [Theorem 1 \cite{Svi2}] \label{IPI3}
Let $F$ be a field of characteristic zero. Then
a non-zero $*$T-ideal of $*$-identities of a finitely generated
associative $F$-algebra with involution coincides
with the $*$T-ideal of $*$-identities of some finite dimensional
associative $F$-algebra with involution.
\end{theorem}
Observe that Theorem \ref{IPI3} can be considered as a partial case of Theorem \ref{IPI2}
assuming that the $\mathbb{Z}_4$-grading is trivial.

\section{Specht problem.}

Theorem \ref{IPI2} yields that any associative
algebra with involution over a field of characteristic zero
has a finite base of $*$-identities.

\begin{theorem} \label{IPI4}
Let $F$ be a field of characteristic zero. Any $*$T-ideal
of the free associative $F$-algebra with involution $F\langle Y, Z \rangle$
is finitely generated as a $*$T-ideal.
\end{theorem}
\noindent {\bf Proof.}
It is clear that $F\langle Y, Z \rangle$ is generated as a $*$T-ideal
by the set $\{ y_1, z_1 \},$ and the zero ideal is generated by the zero polynomial.
Hence it is enough to prove the theorem for proper $*$T-ideals.

Suppose that there exists a proper $*$T-ideal
$\Gamma \subseteq F\langle Y, Z \rangle$ which can not be finitely generated as a
$*$T-ideal.  Then there exists an infinite
sequence of multilinear $*$-polynomials $\{
f_i(x_1,\dots,x_{n_i}) \}_{i \in \mathbb{N}} \subseteq \Gamma,$
such that $\deg f_i < \deg f_j$ for any $i < j,$ and
$f_i \notin *T[f_1,\dots,f_{i-1}]$ for any $i \in \mathbb{N},$ where
$x_j \in Y \cup Z.$

Given $i \in \mathbb{N}$ let us take the $*$T-ideal $\Gamma_i \subseteq F\langle Y, Z
\rangle$ generated by all consequences of the polynomial $f_i$ of
degrees strictly greater then $n_i=\deg f_i.$ Consider
the $*$T-ideal $\widetilde{\Gamma}=\sum_{i \in \mathbb{N}}
\Gamma_i.$ It is clear that for any $i \in \mathbb{N}$ we have that $f_i \notin
\widetilde{\Gamma}.$ By Theorem \ref{IPI2} $\widetilde{\Gamma}$ is the ideal of
identities with involution of the Grassmann $\mathbb{Z}_4$-envelope $E_4(C)$ of some
finite dimensional over $F$ $\mathbb{Z}_4$-graded algebra $C$ with graded involution.

By Lemma 3.1 \cite{Svi3} $C=B \oplus J,$ where $B$ is a
$\mathbb{Z}_4$-graded semisimple algebra with a graded involution,
and $J=J(C)$ is a $\mathbb{Z}_4$-graded nilpotent ideal of $C.$ By
\cite{BahtZaicSeg}, \cite{ZaicSeg} $B$ has the unit $1_B \in
B_{\bar{0}},$ and $1_B$ is symmetric in respect to involution.
Therefore, $E_4(C)=E_4(B) \oplus E_4(J),$ where $E_4(B)$
is a $*$-subalgebra of $E_4(C),$ and $E_4(J)$ is a
nilpotent $*$-ideal of $E_4(C)$ of degree $\mathcal{N}.$

Let us take a polynomial $f_k(x_1,\dots,x_{n_k})$ of degree
$n_k=\deg f_k > \mathcal{N}.$ Consider any evaluation of the polynomial $f_k$
of the type $x_i =a_i=c_{\theta_i} \otimes g_{\theta_i},$ \
$a_i \in E_4(C)^{+}$ if $x_i \in Y,$ and  $a_i \in E_4(C)^{-}$ if $x_i \in Z.$
Where $c_{\theta_i} \in
(B_{\theta_i}^{+} \cup B_{\theta_i}^{-}) \cup (J_{\theta_i}^{+} \cup J_{\theta_i}^{-}),$ \
$g_{\theta_i} \in E_{\theta_i},$ \ $\theta_i \in \mathbb{Z}_4$ for any
$i=1,\dots,n_k.$

If the element $c_{\theta_i}
\in J_{\theta_i}^{+} \cup J_{\theta_i}^{-}$ is radical for any $i=1,\dots,n_k$
then $f_k(a_1,\dots,a_{n_k})=0$ in
$E_4(C),$ since $n_k > \mathcal{N}.$ Suppose that at least one of the
variables of $f_k$ is evaluated by an element of the type $b_{\theta} \otimes
g_{\theta},$ where $b_{\theta} \in B_{\theta}$ is a semisimple element of $C.$
Assume that $x_{\hat{r}}=a_{\hat{r}}=
b_{\theta_{\hat{r}}} \otimes
g_{\theta_{\hat{r}}}$ for any
$b_{\theta_{\hat{r}}} \in B_{\theta_{\hat{r}}},$ \ $g_{\theta_{\hat{r}}}
\in E_{\theta_{\hat{r}}},$ admitting $a_{\hat{r}} \in E_4(C)^{+}$ for $x_{\hat{r}} \in Y$
and $a_{\hat{r}} \in E_4(C)^{-}$ for $x_{\hat{r}} \in Z.$

The algebra $E_4(C)$ has the
natural structure of $E_{\bar{0}}$-module defined by $(c_{\theta} \otimes g_{\theta})
g=c_{\theta}\otimes (g_{\theta} g),$ \ $g \in E_{\bar{0}},$ \ $c_{\theta} \in C_{\theta},$
\ $g_{\theta} \in E_{\theta},$ \ $\theta \in \mathbb{Z}_4.$
This structure preserves the $\mathbb{Z}_4$-grading and the involution
($\deg_{\mathbb{Z}_4} c_{\theta} \otimes (g_{\theta} g) = \deg_{\mathbb{Z}_4} c_{\theta} \otimes g_{\theta} = \theta,$ \ $((c_{\theta} \otimes g_{\theta}) g)^{*} = (c_{\theta} \otimes g_{\theta})^{*} g$). Moreover, $E_{\bar{0}}$ is the center of $E.$
Hence for an element $\tilde{g}_0 \in E_0$ we obtain that
\begin{eqnarray}\label{finSub}
&&f_k(a_1,\dots,a_{\hat{r}},\dots,a_{n_k})
(2 \tilde{g}_0)=f_k(a_1,\dots,b_{\theta_{\hat{r}}}
\otimes (g_{\theta_{\hat{r}}} \cdot
 2 \tilde{g}_0),\dots,a_{n_k})= \\
&&f_k(a_1,\dots,(b_{\theta_{\hat{r}}} \otimes
g_{\theta_{\hat{r}}}) \circ (1_B \otimes
\tilde{g}_0),\dots,a_{n_k})=\tilde{f}_k(a_1,\dots,a_{n_k},1_B
\otimes \tilde{g}_0)=0. \nonumber
\end{eqnarray}
Here $\tilde{f}_k(x_1,\dots,x_{n_k},y_{\bar{0}})=$ $f_k(x_1,\dots,x_{\hat{r}} \circ
y_{\bar{0}},\dots,x_{n_k}) \in
\widetilde{\Gamma}=\mathrm{Id}^{*}(E_4(C)),$ $y_{\bar{0}} \in
Y_{\bar{0}}.$
The equality $f_k(a_1,\dots,a_{n_k})=0$ directly follows from (\ref{finSub}).

Since $f_k$ is multilinear then it implies that $f_k \in
\mathrm{Id}^{*}(E_4(C))=\widetilde{\Gamma}.$ This
contradicts to the construction of $\widetilde{\Gamma}.$
Therefore, $\Gamma$ is finitely generated as a $*$T-ideal.
\hfill $\Box$

Observe that the usual Grassmann envelope of the superalgebras with
superinvolution also can be considered in the context of the Specht problem
and Classification theorems for identities with involution.
We assume that results similar to Theorem 6.2 \cite{Svi3}, and
Theorems \ref{IPI1}, \ref{IPI2} can be obtained also in this case.

\begin{conjecture} \label{1}
Let $F$ be a field of characteristic zero, and
$A=A_{\bar{0}} \oplus A_{\bar{1}}$ a finitely generated associative PI-superalgebra
over $F$ with superinvolution. Then there exists a finite dimensional over $F$
associative superalgebra $C=C_{\bar{0}} \oplus C_{\bar{1}}$ with superinvolution
which satisfies the same identities with superinvolution as $A.$
\end{conjecture}

\begin{conjecture} \label{2}
Let $F$ be a field of characteristic zero. Then any associative $F$-algebra with
involution satisfies the same $*$-identities as the Grassmann envelope
$E(C)=C_{\bar{0}} \otimes E_{\bar{0}} \oplus C_{\bar{1}} \otimes E_{\bar{1}}$ of
some associative superalgebra $C=C_{\bar{0}} \oplus C_{\bar{1}}$ with superinvolution,
finite dimensional over $F.$
\end{conjecture}

The confirmation of these conjectures could imply another solution of the Specht
problem for $*$-identities.

Author is deeply thankful to Ivan Shestakov and Antonio Giambruno
for useful discussions and inspiration, and grateful to FAPESP for the
financial support in this work.

\bibliographystyle{amsplain}

\begin{thebibliography}{99}
\bibitem {AB} E.Aljadeff, A.Kanel-Belov, Representability and
Specht problem for $G$-graded algebras, Adv. math., 225(2010), 2391-2428.

\bibitem {A3} S.A.Amitsur, Rings with involution, Israel J. Math., 6(1968), 99-106.

\bibitem {A4} S.A.Amitsur, Identities in rings with involution, Israel J. Math., 7(1968), 63-68.

\bibitem{BahtGR} Yu.Bakhturin, A.Giambruno, D.Riley, Group-graded
algebras with polynomial identities, Israel J. Math., 104(1998), 145-155.

\bibitem{BahtZaicSeg} Yu.A.Bakhturin, S.K.Sehgal, M.V.Zaicev,
Finite-dimensional simple graded algebras, Sb. Math., 199(7)(2008), 965-983.

\bibitem {BahtShestZ} Y.A.Bakhturin, I.P.Shestakov, M.V.Zaicev, Gradings
on simple Jordan and Lie algebras, J. Algebra, 283(2005), 849-868.

\bibitem {B1} A.Ya.Belov, Local finite basis property and local representability of varieties
of associative rings, (Russian), Izv. Ross. Akad. Nauk, Ser. Mat., 74(2010), no. 1, 3-134;
English transl. in Izv. Math., 74(2010), no. 1, 1-126.

\bibitem {B2} A.Ya.Belov, Local finite basis property and local finite representability of
varieties of associative rings. (Russian) Dokl. Akad. Nauk, 432(2010), no. 6, 727-731;
English transl. in Dokl. Math. 81(2010), no. 3, 458-461.

\bibitem {B4} A.Belov-Kanel, L.Rowen, U.Vishne, Structure of Zariski-closed algebras,
Trans. Amer. Math. Soc., 362(2010), no. 9, 4695-4734.

\bibitem {B5} A.Belov, L.H.Rowen, U.Vishne, Application of full quivers of representations
of algebras to polynomial identities, Comm. Algebra, 39(2011), no. 12, 4536-4551.

\bibitem {B6} A.Belov, L.H.Rowen, U.Vishne, Full exposition of Specht's problem,
Serdica Math. J., 38(2012), no. 1-3, 313-370.

\bibitem {B7} A.Belov, L.H.Rowen, U.Vishne, Full quivers of representations of algebras,
Trans. Amer. Math. Soc., 364(2012), no. 10, 5525-5569.

\bibitem {B8} A.Belov, L.H.Rowen, U.Vishne, PI-varieties associated to full
quivers of representations of algebras, Trans. Amer. Math. Soc., 365(2013), no. 5, 2681-2722.

\bibitem {B9} A.Kanel-Belov, L.H.Rowen, U.Vishne, Specht's problem for associative
affine algebras over commutative Noetherian rings, Trans. Amer. Math. Soc., in press.

\bibitem{BergC} J.Bergen, M.Cohen, Action of commutative Hopf
algebras, Bull. London Math. Soc., 18(1986), 159-164.

\bibitem{CurtR} C.W.Curtis, I.Reiner, Representation theory of finite groups and associative algebras,
Reprint of the 1962 original, AMS Chelsea Publishing, Providence, RI, 2006.

\bibitem {Dren} V.Drensky, Free algebras and PI-algebras, Springer-Verlag Singapore,
Singapore, 2000.

\bibitem {DrenForm} V.Drensky, E.Formanek, Polynomial identity rings,
Birkhauser Verlag, Basel-Boston-Berlin, 2004.

\bibitem {DrenGiam} V.Drensky, A.Giambruno, Cocharacters, codimensions and Hilbert series
of the polynomial identities for $2 \times 2$ matrices with involution,
Canad. J. Math., 46(1994), 718–733.

\bibitem {GM} A.Giambruno, S.Mishchenko, On star-varieties with almost
polynomial growth, Algebra Colloquium, 8(1)(2001), 33-42.

\bibitem {GM1} A.Giambruno, S.Mishchenko, Super-cocharacters, star-cocharacters
and multiplicities bounded by one, Manuscripta Math. 128(2009), 483–504.

\bibitem {GR} A.Giambruno, A.Regev, Wreath products and P.I. algebras,
J. Pure Applied Algebra, 35(1985), 133-149.

\bibitem{GRZ1} A.Giambruno, A.Regev, M.Zaicev, Polynomial idintities
and Combinatorial Methods, Marcel Dekker Inc., New York, Basel,
2003.

\bibitem {GZbook} A.Giambruno, M.Zaicev, Polynomial identities and
asymptotic methods, Amer.Math.Soc., Math. Surveys and
Monographs, 122, Providence, R.I., 2005.

\bibitem{LieHum} J.Humphreys, Introduction to Lie Algebras and Representation Theory,
Grad. Textbooks, second ed., Springer-Verlag, Berlin, 2003.

\bibitem{JKer} G.James, A.Kerber, The Representation Theory of the Symmetric Group,
Encyclopedia of Mathematics and Its Applications 16, Addison-Wesley, London, 1981.

\bibitem {BelRow} A.Kanel-Belov, L.H.Rowen, Computational aspects
of polynomial identities, A K Peters Ltd., Wellesley, MA, 2005.

\bibitem {Kem1} A.R.Kemer, Ideals of identities of associative
algebras, Amer.Math.Soc. Translations of Math. Monographs,
87, Providence, R.I., 1991.

\bibitem{Invbook} M.-A.Knus, A.A.Merkurjev, M.Rost, J.-P.Tignol,
The Book of Involutions, Colloquium Publications, 44, AMS, 1998.

\bibitem{McCrim} K.McCrimmon, A taste of Jordan algebras, Universitext,
Berlin, New York, Springer-Verlag, 2004.

\bibitem{Reg1} A.Regev, Existence of identities in $A\otimes B$, Israel
J. Math., 11(1972), 131-152.

\bibitem{Specht} W.Specht, Gesetze in Ringen I, Math. Z., 52(5)(1950), 557-589.

\bibitem {Svi1} I.Sviridova, Identities of PI-algebras graded
by a finite abelian group, Comm. Algebra, 39(9)(2011),
3462-3490.

\bibitem {Svi2} I.Sviridova, Finitely generated algebras with involution
and their identities, J. Algebra, 383(2013), 144-167.

\bibitem {Svi3} I.Sviridova, Identities of finitely generated
graded algebras with involution, arXiv:1410.2222 [math.RA], preprint.

\bibitem{ZaicSeg} M.V.Zaicev, S.K.Sehgal, Finite gradings on simple
Artinian rings, Vestnik Mosk. Univ., Matem., Mechan., 3(2001),
21-24 (Russian).

\bibitem{ZSSS} K.A.Zhevlakov, A.M.Slin'ko, I.P.Shestakov, A.I.Shirshov, Rings
that are nearly associative, Pure and Applied Mathematics, 104, Academic Press, Inc. New
York-London, 1982.
\end{thebibliography}

\end{document}